\documentclass[12pt]{amsart}

\usepackage{enumitem}
\usepackage{psfrag}
\usepackage{graphicx}
\usepackage{pinlabel}
\usepackage{graphicx}
\usepackage{amsmath}
\usepackage{amssymb}
\usepackage{amscd}
\usepackage{autobreak}
\usepackage{picinpar}
\usepackage{times}
\usepackage{pb-diagram}
\usepackage{graphicx}
\usepackage{wrapfig}
\usepackage{xspace}
\usepackage{hyperref}
\usepackage{url}
\usepackage{caption}
\usepackage{subcaption}
\usepackage{fancyhdr}
\usepackage{overpic}
\usepackage{color}
\usepackage[square,comma,sort&compress,numbers]{natbib}
\usepackage{latexsym}
\usepackage{mathrsfs}
\usepackage{amsfonts}
\usepackage{hyperref}
\usepackage{amsthm}
\usepackage{appendix}
\usepackage{pdfsync}

\theoremstyle{definition}
\newtheorem{theorem}{Theorem}[section]
\newtheorem{lemma}[theorem]{Lemma}
\newtheorem{proposition}[theorem]{Proposition}
\newtheorem{corollary}[theorem]{Corollary}
\newtheorem{definition}[theorem]{Definition}

\newtheorem{remark}[theorem]{Remark}
\newtheorem*{theorem*}{Theorem}
\def\qed{\hfill{Q.E.D.}\smallskip}

\begin{document}

\title{\bf A discrete uniformization theorem for decorated piecewise Euclidean metrics on surfaces}
\author{Xu Xu, Chao Zheng}

\date{\today}

\address{School of Mathematics and Statistics, Wuhan University, Wuhan, 430072, P.R.China}
 \email{xuxu2@whu.edu.cn}

\address{School of Mathematics and Statistics, Wuhan University, Wuhan 430072, P.R. China}
\email{czheng@whu.edu.cn}



\begin{abstract}
In this paper, we introduce a new discretization of the Gaussian curvature on surfaces,
which is defined as the quotient of the angle defect and the area of  some dual cell of a weighted triangulation at the conic singularity.
A discrete uniformization theorem for this discrete Gaussian curvature is established on surfaces with non-positive Euler number.
The main tools are Bobenko-Lutz's discrete conformal theory for decorated piecewise Euclidean metrics on surfaces and variational principles with constraints.
\end{abstract}

\maketitle

\section{Introduction}

Bobenko-Lutz \cite{BL} recently introduced the decorated piecewise Euclidean metrics on surfaces.
Suppose $S$ is a connected closed surface and $V$ is a finite non-empty subset of $S$, we call $(S, V)$ a marked surface.
A piecewise Euclidean metric (PE metric) $dist_{S}$ on the marked surface $(S,V)$ is a flat cone metric
with the conic singularities contained in $V$.
A decoration on a PE surface $(S,V, dist_{S})$ is a choice of circle of radius $r_i\geq 0$ at each point $i\in V$.
These circles in the decoration are called vertex-circles.
We denote a decorated PE surface by $(S,V, dist_{S}, r)$ and call the pair $(dist_S,r)$ a decorated PE metric.
In this paper, we focus on the case that $r_i>0$ for all $i\in V$ and each pair of vertex-circles is separated.

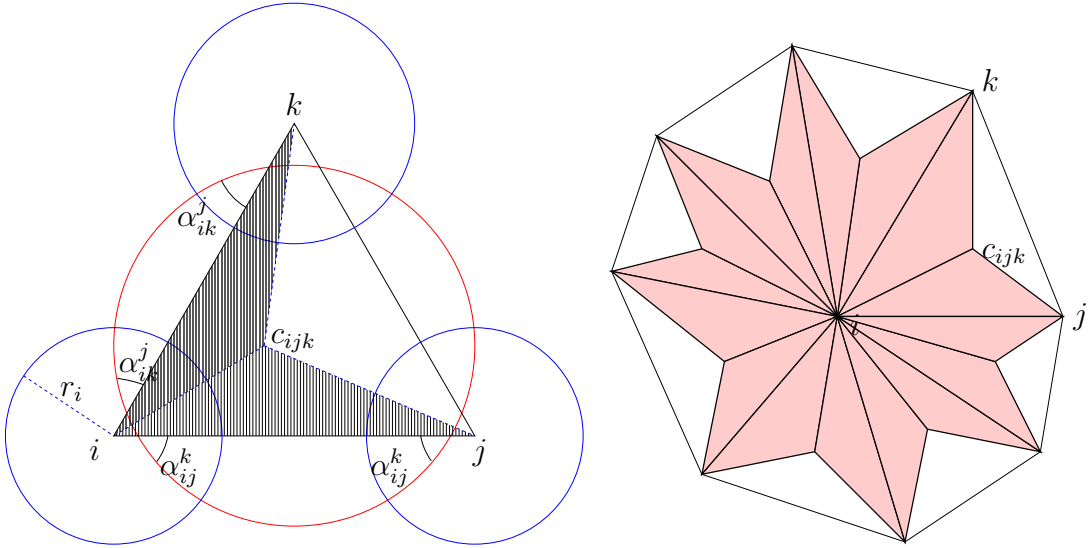
\begin{figure}[htbp]
\centering
\begin{minipage}[t]{0.5\columnwidth}
\begin{overpic}[scale=0.4]{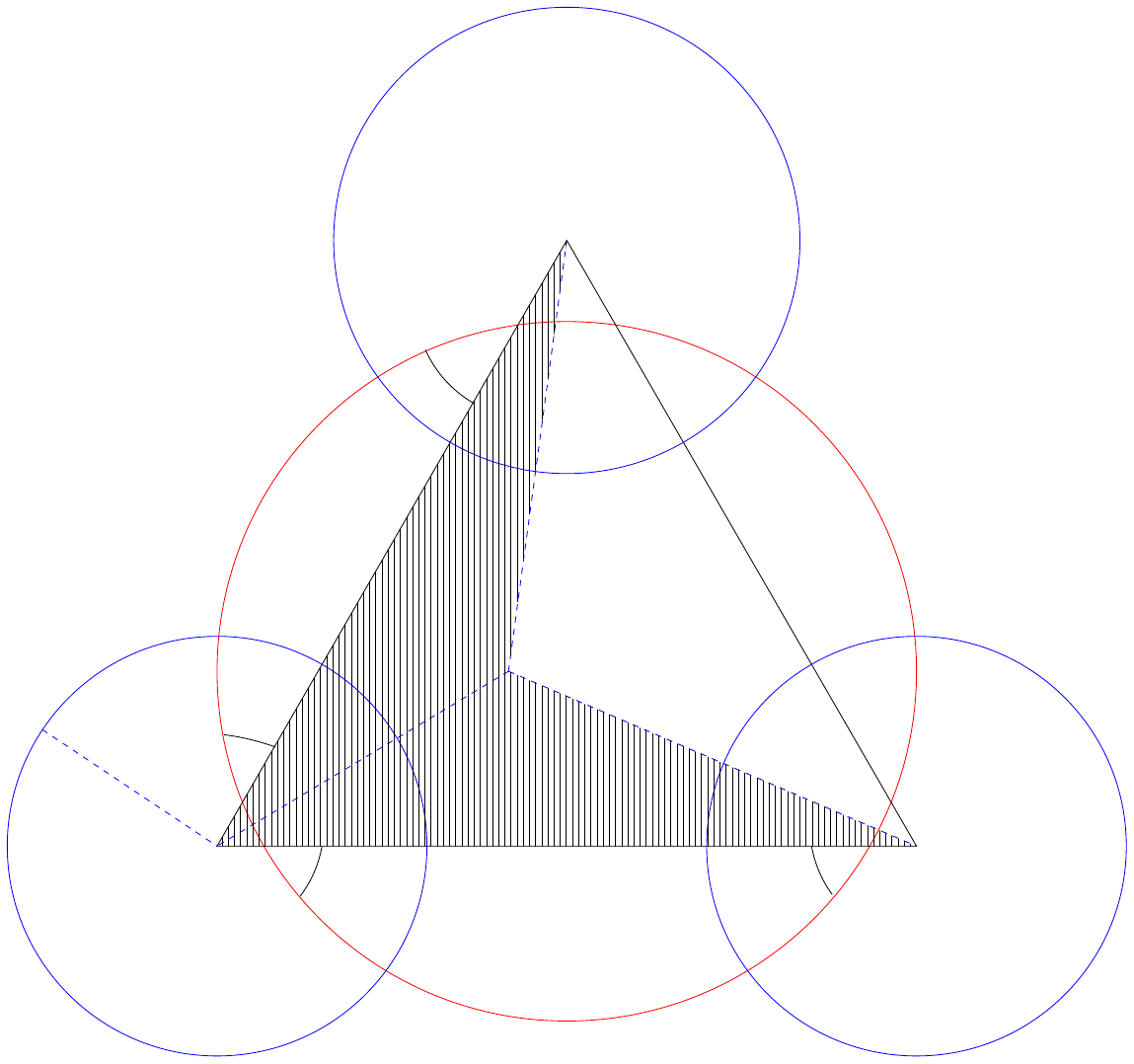}
\put (15,15) {$i$}
\put (80.5,15) {$j$}
\put (48.5,74) {$k$}
\put (10,26) {$r_i$}
\put (46,35) {$c_{ijk}$}
\put (30,55) {$\alpha_{ik}^j$}
\put (20,30) {$\alpha_{ik}^j$}
\put (27,12.5) {$\alpha_{ij}^k$}
\put (63,12.5) {$\alpha_{ij}^k$}
\end{overpic}
\end{minipage}
\begin{minipage}[t]{0.4\columnwidth}
\begin{overpic}[scale=0.6]{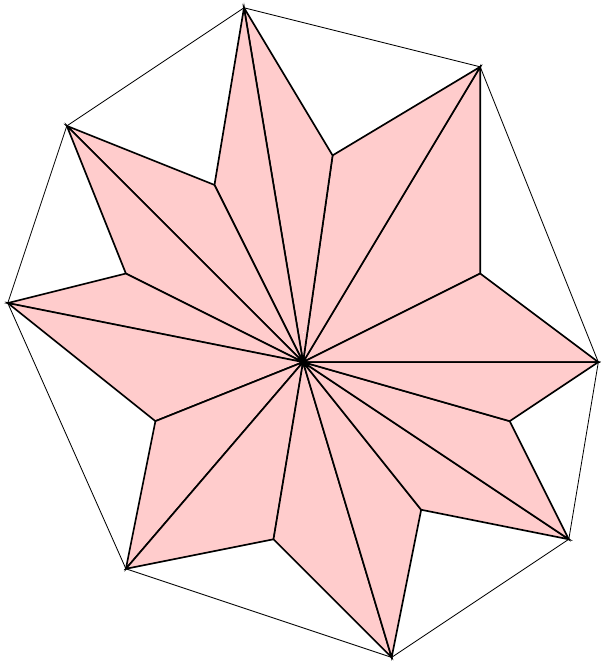}
\put (48,42) {$i$}
\put (92,44) {$j$}
\put (74,90) {$k$}
\put (74,57) {$c_{ijk}$}
\end{overpic}
\end{minipage}
\caption{Domain of the signed area $A_i^{jk}$ in a decorated triangle $\{ijk\}$ (left) and domain of the area $A_i$ in a decorated PE surface (right)}
\label{figure 1}
\end{figure}

For a PE surface $(S,V, dist_{S})$, denote $\theta_i$ as the cone angle at $i\in V$.
The angle defect
\begin{equation}\label{Eq: curvature W}
W: V\rightarrow (-\infty, 2\pi), \quad W_i=2\pi-\theta_i,
\end{equation}
is used to describe the conic singularities of the PE metric.
Let $\mathcal{T}={(V,E,F)}$ be a triangulation of $(S, V)$, where $V,E,F$ are the sets of vertices, edges and faces respectively.
The triangulation $\mathcal{T}$ for a PE surface $(S,V, dist_S)$ is a geodesic triangulation if the edges are geodesics in the PE metric $dist_S$.
We use one index to denote a vertex (such as $i$), two indices to denote an edge (such as $\{ij\}$) and three indices to denote a face (such as $\{ijk\}$) in the triangulation $\mathcal{T}$.
For any decorated geodesic triangle $\{ijk\}\in F$,
there is a unique circle $C_{ijk}$ simultaneously orthogonal to the three vertex-circles at the vertices $i,j,k$ \cite{Glickenstein preprint}.
We call this circle $C_{ijk}$ as the face-circle of the decorated geodesic triangle $\{ijk\}$ and denote its center by $c_{ijk}$ and radius by $r_{ijk}$.
The center $c_{ijk}$ of the face-circle $C_{ijk}$ of the decorated geodesic triangle $\{ijk\}$ is the geometric center introduced by Glickenstein \cite{Glickenstein JDG} and Glickenstein-Thomas \cite{GT}
for general discrete conformal structures on surfaces.
Denote $\alpha_{ij}^k$ as the interior intersection angle of the face-circle $C_{ijk}$ and the edge $\{ij\}$.
Please refer to Figure \ref{figure 1} (left) for the angle $\alpha_{ij}^k$.
The edge $\{ij\}$, shared by two adjacent decorated triangles $\{ijk\}$ and $\{ijl\}$, is called weighted Delaunay if
\begin{equation}\label{Eq: F7}
\alpha_{ij}^k+\alpha_{ij}^l\leq \pi.
\end{equation}
The triangulation $\mathcal{T}$ is called weighted Delaunay in the decorated PE metric $(dist_S,r)$ if every edge in the triangulation is weighted Delaunay.
Connecting the center $c_{ijk}$ with the vertices $i,j,k$ by geodesics produces a cellular decomposition of the decorated triangle $\{ijk\}$.
Denote $A_i^{jk}$ as the sum of the signed area of the two triangles adjacent to $i$ in the cellular decomposition of the decorated triangle $\{ijk\}$.
The area of the triangle with the vertices $i$, $j$, $c_{ijk}$ is positive if it is on the same side of the edge $\{ij\}$ as the decorated  triangle $\{ijk\}$,
otherwise it is negative (or zero if $c_{ijk}$ lies in $\{ij\}$).
Please refer to the shaded domain in Figure \ref{figure 1} (left) for $A_i^{jk}$.
Gluing these cells of all decorated triangles isometrically along edges in pairs leads a cellular decomposition of the decorated PE surface $(S,V, dist_S,r)$.
Set
$$A_i=\sum_{\{ijk\}\in F}A_i^{jk}.$$
Please refer to Figure \ref{figure 1} (right) for $A_i$.

\begin{definition}\label{Def: new curvature}
Suppose  $(S,V, dist_S, r)$ is a decorated PE surface and  $\mathcal{T}$ is a weighted Delaunay triangulation of $(S,V, dist_S, r)$.
The discrete Gaussian curvature $K_i$ at the vertex $i\in V$ is the quotient of the angle defect $W_i$ and the area $A_i$ of the dual cell at the vertex $i\in V$, i.e.,
\begin{equation}\label{Eq: K_i}
K_i=\frac{W_i}{A_i}.
\end{equation}
\end{definition}

\begin{remark}
In the literature, the discrete curvature is usually defined by the angle defect $W$ in (\ref{Eq: curvature W}).
However, the angle defect $W$ is scaling invariant and does not approximate the smooth Gaussian curvature pointwisely on smooth surfaces as the triangulations of the surface become finer and finer.
This is supported by the discussions in \cite{BPS, GX2}.
For the discrete Gaussian curvature $K$ in (\ref{Eq: K_i}), it scales by a factor $\frac{1}{u^2}$ upon a global rescaling of the decorated PE metric by a factor $u$.
This property is paralleling to that of the smooth Gaussian curvature on surfaces.
On the other hand, the definition of the discrete Gaussian curvature $K_i$ coincides with the original definition of the Gaussian curvature on smooth surfaces.
This implies that the discrete Gaussian curvature $K_i$
is a good candidate as a discretization of the smooth Gaussian curvature on surfaces.
\end{remark}

\begin{remark}
According to Definition \ref{Def: new curvature},
the discrete Gaussian curvature $K_i$ defined by (\ref{Eq: K_i}) seems to depend on the choice of weighted Delaunay triangulations of
the decorated PE surface $(S, V, dist_S, r)$.
We will show that $K_i$ is an intrinsic geometric invariant of the decorated PE surface $(S,V, dist_S, r)$
in the sense that it is independent of the weighted Delaunay triangulations of $(S, V, dist_S, r)$.
Note that the angle defect $W_i$ defined by (\ref{Eq: curvature W}) is an intrinsic geometric invariant of a decorated PE surface,
we just need to prove that $A_i$ is independent of the choice of weighted Delaunay triangulations.
This is true by Lemma \ref{Lem: independent}.
\end{remark}

\begin{remark}
The weighted Delaunay triangulation is a natural generalization of the classical Delaunay triangulation.
When the weighted Delaunay triangulation is reduced to the classical Delaunay triangulation, i.e. $r_i=0$ for all $i\in V$,
the area $A_i$ is exactly twice the area of the Voronoi cell at the vertex $i$.
Thus the area $A_i$ is a generalization of the area of the Voronoi cell at the vertex $i$.
As a result, the discrete Gaussian curvature in Definition \ref{Def: new curvature} generalizes Kou\v{r}imsk\'{a}'s definition of discrete Gaussian curvature in \cite{Kourimska,Kourimska Thesis}.
\end{remark}

The discrete Yamabe problem for a decorated PE metric $(dist_S,r)$ on $(S,V)$ asks if there exists a discrete conformal equivalent decorated PE metric on $(S,V)$ with constant discrete Gaussian curvature.
The following discrete uniformization theorem solves this problem affirmatively for the discrete Gaussian curvature $K$ in Definition \ref{Def: new curvature}.

\begin{theorem}\label{Thm: main}
For any decorated PE metric $(dist_S,r)$ on a marked surface $(S,V)$ with Euler number $\chi(S)\leq 0$,
there is a discrete conformal equivalent decorated PE metric with constant discrete Gaussian curvature $K$.
\end{theorem}

By the relationships of the discrete Gaussian curvature $K$ and the classical discrete Gaussian curvature $W$,
the case $\chi(S)=0$ in Theorem \ref{Thm: main} is covered by Bobenko-Lutz's work \cite{BL}.
Therefore, we just need to prove the case $\chi(S)<0$ in Theorem \ref{Thm: main}.

\begin{remark}
The discrete Yamabe problem on surfaces for different types of discrete conformal structures with respect to the classical discrete Gaussian curvature $W$ has been extensively studied in the literature.
For Thurston's circle packings on surfaces,
the solution of discrete Yamabe problem gives rise to the famous Koebe-Andreev-Thurston Theorem.
See also the work of Beardon-Stephenson \cite{BS uniformization} for the discrete uniformization theorems for circle packings on surfaces.
For the vertex scalings introduced by Luo \cite{Luo1} on surfaces,
Gu-Luo-Sun-Wu \cite{Gu1}, Gu-Guo-Luo-Sun-Wu \cite{Gu2}, Springborn \cite{Springborn} and Izmestiev-Prosanov-Wu \cite{IPW} give nice answers to this problem in different background geometries.
Recently, Bobenko-Lutz \cite{BL} established the discrete conformal theory for decorated PE metrics and prove the corresponding discrete uniformization theorem.
Since Bobenko-Lutz's discrete conformal theory of decorated PE metrics also applies to the Euclidean vertex scalings and thus generalizes Gu-Luo-Sun-Wu's result \cite{Gu1} and Springborn's result \cite{Springborn},
Theorem \ref{Thm: main} also generalizes Kou\v{r}imsk\'{a}'s results in \cite{Kourimska,Kourimska Thesis}.
It should be mentioned that Kou\v{r}imsk\'{a} \cite{Kourimska,Kourimska Thesis} constructed counterexamples to the uniqueness of PE metrics with constant discrete Gaussian curvatures.
We conjecture that the decorated PE metric with constant discrete Gaussian curvature $K$ in Theorem \ref{Thm: main} is not unique.
\end{remark}

The main tools for the proof of Theorem \ref{Thm: main} are
Bobenko-Lutz's discrete conformal theory for decorated PE metrics on surfaces \cite{BL} and variational principles with constraints.
The main ideas of the paper come from reading of Bobenko-Lutz \cite{BL} and Kou\v{r}imsk\'{a} \cite{Kourimska,Kourimska Thesis}.

The paper is organized as follows.
In Section \ref{Sec: Preliminaries},
we briefly recall Bobenko-Lutz's discrete conformal theory for decorated PE metrics on surfaces.
Then we show that $A_i$ is independent of the choice of weighted Delaunay triangulations, i.e., Lemma \ref{Lem: independent}.
We also give some notations and  a variational characterization of the area $A_i^{jk}$.
In this section, we also extend the energy function $\mathcal{E}$ and the area function $A_{tot}$.
In Section \ref{Sec: proof the main theorem},
we translate Theorem \ref{Thm: main} into an optimization problem with constraints, i.e., Lemma \ref{Lem: minimum lies at the boundary}.
Using the classical result from calculus, i.e., Theorem \ref{Thm: calculus}, we translate Lemma \ref{Lem: minimum lies at the boundary} into Theorem \ref{Thm: key}.
By analysing the limit behaviour of sequences of discrete conformal factors, we get an asymptotic expression of the function $\mathcal{E}$, i.e., Lemma \ref{Lem: E decomposition}.
In the end, we prove Theorem \ref{Thm: key}.
\\
\\
\textbf{Acknowledgements}\\[8pt]
The first author thanks Professor Feng Luo for his invitation to the workshop
``Discrete and Computational Geometry, Shape Analysis, and Applications" taking place
at Rutgers University, New Brunswick from May 19th to May 21st, 2023.
The first author also thanks Carl O. R. Lutz for helpful communications during the workshop.

\section{Preliminaries on decorated PE surfaces}\label{Sec: Preliminaries}

\subsection{Discrete conformal equivalence and Bobenko-Lutz's discrete conformal theory}
In this subsection, we briefly recall Bobenko-Lutz's discrete conformal theory for decorated PE metrics on surfaces.
Please refer to Bobenko-Lutz's original work \cite{BL} for more details on this.
The PE metric $dist_{S}$ on a PE surface with a geodesic triangulation defines a length map $l: E\rightarrow \mathbb{R}_{>0}$ such that $l_{ij}, l_{ik}, l_{jk}$ satisfy the triangle inequalities for any triangle $\{ijk\}\in F$.
Conversely, given a function $l: E\rightarrow \mathbb{R}_{>0}$ satisfying the triangle inequalities for any face $\{ijk\}\in F$,
one can construct a PE metric on a triangulated surface by isometrically gluing Euclidean triangles along edges in pairs.
Therefore, we use $l: E\rightarrow \mathbb{R}_{>0}$ to denote a PE metric and use $(l,r)$ to denote a decorated PE metric on a triangulated surface $(S,V,\mathcal{T})$.

\begin{definition}[\cite{BL}, Proposition 2.2]
\label{Def: DCE}
Let $\mathcal{T}$ be a triangulation of a marked surface $(S,V)$.
Two decorated PE metrics $(l,r)$ and $(\widetilde{l},\widetilde{r})$ on $(S,V, \mathcal{T})$ are discrete conformal equivalent
if and only if there exists a discrete conformal factor $u\in \mathbb{R}^V$ such that
\begin{equation}\label{Eq: DCE1}
\widetilde{r}_i=e^{u_i}r_i,
\end{equation}
\begin{equation}\label{Eq: DCE2}
\widetilde{l}_{ij}^2
=(e^{2u_i}-e^{u_i+u_j})r^2_i
+(e^{2u_j}-e^{u_i+u_j})r^2_j
+e^{u_i+u_j}l_{ij}^2
\end{equation}
for all $\{ij\}\in E$.
\end{definition}

\begin{remark}
Note that the inversive distance
\begin{equation}\label{Eq: inversive distance}
I_{ij}=\frac{l^2_{ij}-r^2_i-r^2_j}{2r_ir_j}
\end{equation}
between two vertex-circles is invariant under M\"{o}bius transformations \cite{Coxeter}.
Combining (\ref{Eq: DCE1}) and (\ref{Eq: DCE2}) gives $I=\widetilde{I}$.
Since each pair of vertex-circles is required to be separated, we have $I>1$.
Therefore, Definition \ref{Def: DCE} can be regarded as a special case of the inversive distance circle packings introduced by Bowers-Stephenson \cite{BS}.
One can refer to \cite{CLXZ, Guo,Luo GT,Xu AIM,Xu MRL} for more properties of the inversive distance circle packings on triangulated surfaces.
\end{remark}

In general, the existence of decorated PE metrics with constant discrete Gaussian curvatures on triangulated surfaces can not be guaranteed if the triangulation is fixed.
In the following, we work with a generalization of the discrete conformal equivalence in Definition \ref{Def: DCE},
introduced by Bobenko-Lutz \cite{BL},
which allows the triangulation of the marked surface to be changed under the weighted Delaunay condition.

\begin{definition}[\cite{BL}, Definition 4.11]
\label{Def: GDCE}
Two decorated PE metrics $(dist_{S},r)$ and $(\widetilde{dist}_{S},\widetilde{r})$ on the marked surface $(S,V)$ are discrete conformal equivalent if there is a sequence of triangulated decorated PE surfaces
$(\mathcal{T}^0,l^0,r^0),...,(\mathcal{T}^N,l^N,r^N)$ such that
\begin{description}
  \item[(i)] the decorated PE metric of $(\mathcal{T}^0,l^0,r^0)$ is $(dist_{S},r)$ and the decorated PE metric of $(\mathcal{T}^N,l^N,r^N)$ is $(\widetilde{dist}_{S},\widetilde{r})$,
  \item[(ii)] each $\mathcal{T}^n$ is a weighted Delaunay triangulation of the decorated PE surface $(\mathcal{T}^n,l^n,r^n)$,
  \item[(iii)] if $\mathcal{T}^n=\mathcal{T}^{n+1}$, then there is a discrete conformal factor $u\in \mathbb{R}^V$ such that $(\mathcal{T}^n,l^n,r^n)$ and $(\mathcal{T}^{n+1},l^{n+1},r^{n+1})$ are related by (\ref{Eq: DCE1}) and (\ref{Eq: DCE2}),
  \item[(iv)] if $\mathcal{T}^n\neq\mathcal{T}^{n+1}$, then $\mathcal{T}^n$ and $\mathcal{T}^{n+1}$ are two different weighted Delaunay triangulations of the same decorated PE surface.
\end{description}
\end{definition}

Definition \ref{Def: GDCE} defines an equivalence relationship for decorated PE metrics on a marked surface.
The equivalence class of a decorated PE metric $(dist_S,r)$ on $(S,V)$ is also called as the discrete conformal class of $(dist_S,r)$ and denoted by $\mathcal{D}(dist_S,r)$.

\begin{lemma}[\cite{BL}]
The discrete conformal class $\mathcal{D}(dist_S,r)$ of a decorated PE metric $(dist_S,r)$ on the marked surface $(S,V)$ is parameterized by $\mathbb{R}^V=\{u: V\rightarrow \mathbb{R}\}$.
\end{lemma}
For simplicity, for any $(\widetilde{dist}_S,\widetilde{r})\in \mathcal{D}(dist_S,r)$,
we denote it by $(dist_S(u),r(u))$ for some $u\in \mathbb{R}^V$.
Set
\begin{equation*}
\mathcal{C}_\mathcal{T}(dist_{S},r)
=\{u\in \mathbb{R}^V |\ \mathcal{T}\ \text{is a weighted Delaunay triangulation of}\ (S,V,dist_S(u),r(u))\}.
\end{equation*}

For any decorated PE surface, there exists a unique complete hyperbolic surface $\Sigma_g$, i.e., the hyperbolic surface induced by any triangular refinement of its unique weighted Delaunay tessellation.
It is homeomorphic to $S\backslash V$ and called as the fundamental discrete conformal invariant of the decorated PE metric $(dist_{S},r)$.
The decoration of $\Sigma_g$ is denoted by $\omega:=e^{h}$ and here the height $h$ is related to $u$ by $dh_i=-du_i$.
The canonical weighted Delaunay tessellation $\mathcal{T}$ of $\Sigma_g$ is denoted by $\mathcal{T}_{\Sigma_g}^\omega$.
Bobenko-Lutz \cite{BL} defined the following set
\begin{equation*}
\mathcal{D}_\mathcal{T}(\Sigma_g)
=\{\omega\in \mathbb{R}_{>0}^V|\mathcal{T}\ \text{refines}\ \mathcal{T}_{\Sigma_g}^\omega\}
\end{equation*}
and proved the following proposition.

\begin{proposition}[\cite{BL}, Proposition 4.3]
\label{Prop: finite decomposition}
Given a complete hyperbolic surface with ends $\Sigma_g$.
\begin{description}
\item[(1)]
Each $\mathcal{D}_{\mathcal{T}_n}(\Sigma_g)$ is either empty or the intersection of $\mathbb{R}^V_{>0}$ with a closed polyhedral cone.
\item[(2)]
There is only a finite number of geodesic tessellations $\mathcal{T}_1,...,\mathcal{T}_N$ of $\Sigma_g$ such that $\mathcal{D}_{\mathcal{T}_n}(\Sigma_g)$ $(n=1,...,N)$ is non-empty.
In particular, $\mathbb{R}^V_{>0}
=\bigcup_{n=1}^N\mathcal{D}_{\mathcal{T}_n}(\Sigma_g)$.
\end{description}
\end{proposition}

Let $P$ be the polyhedral cusp corresponding to the triangulated surface $(S,V,\mathcal{T})$ with fundamental discrete conformal invariant $\Sigma_g$.
The polyhedral cusp is convex if and only if $\mathcal{T}$ is a weighted Delaunay triangulation.
The set of all heights $h$ of convex polyhedral cusps over the triangulated hyperbolic surface $(\Sigma_g,\mathcal{T})$ is denoted by $\mathcal{P}_\mathcal{T}(\Sigma_g)\subseteq \mathbb{R}^V$.

\begin{proposition}[\cite{BL}, Proposition 4.9]
\label{Prop: some spaces}
Given a decorated PE metric $(dist_{S},r)$ on the marked surface $(S,V)$.
Then $\mathcal{C}_\mathcal{T}(dist_{S},r)$, $\mathcal{P}_\mathcal{T}(\Sigma_g)$ and $\mathcal{D}_\mathcal{T}(\Sigma_g)$ are homeomorphic.
\end{proposition}

Combining Proposition \ref{Prop: finite decomposition} and Proposition \ref{Prop: some spaces} gives the following result.
\begin{lemma}[\cite{BL}]\label{Lem: finite decomposition}
The set
\begin{equation*}
J=\{\mathcal{T}| \mathcal{C}_{\mathcal{T}}(dist_{S},r)\ \text{has non-empty interior in}\ \mathbb{R}^V\}
\end{equation*}
 is a finite set, $\mathbb{R}^V=\cup_{\mathcal{T}_i\in J}\mathcal{C}_{\mathcal{T}_i}(dist_{S},r)$, and
each $\mathcal{C}_{\mathcal{T}_i}(dist_{S},r)$ is homeomorphic to a polyhedral cone (with its apex removed)
and its interior is homeomorphic to $\mathbb{R}^V$.
\end{lemma}

\subsection{A decorated triangle}

Denote $r_{ij}$ as half of the distance of the two intersection points of the face-circle $C_{ijk}$ and the edge $\{ij\}$.
Denote $h_{ij}^k$ as the signed distance of the center $c_{ijk}$ to the edge $\{ij\}$,
which is defined to be positive if the center is on the same side of the line determined by $\{ij\}$ as the triangle $\{ijk\}$ and negative otherwise (or zero if the center is on the line).
Note that $h_{ij}^k$ is symmetric in the indices $i$ and $j$.
By Figure \ref{figure 2}, we have
\begin{equation}\label{Eq: F6}
h_{ij}^k=r_{ij}\cot\alpha_{ij}^k.
\end{equation}
Since $r_{ij}>0$ and $\alpha_{ij}^k\in (0,\pi)$,
if $h_{ij}^k<0$, then $\alpha_{ij}^k\in (\frac{\pi}{2},\pi)$.
The equality (\ref{Eq: F6}) implies that (\ref{Eq: F7}) is equivalent to
\begin{equation}\label{Eq: F9}
h_{ij}^k+h_{ij}^l\geq 0
\end{equation}
for any adjacent triangles $\{ijk\}$ and $\{ijl\}$ sharing a common edge $\{ij\}$.
Therefore, the equality (\ref{Eq: F9}) also characterizes
a weighted Delaunay triangulation $\mathcal{T}$ for a decorated PE metric $(l,r)$ on $(S,V)$.
Due to this fact, the equality (\ref{Eq: F9}) is usually used to define the weighted Delaunay triangulations of decorated PE surfaces.
See \cite{CLXZ,Glickenstein DCG} and others for example.
Then $A_i^{jk}$ can be written as
\begin{equation}\label{Eq: Area 1}
A_i^{jk}=\frac{1}{2}l_{ij}h_{ij}^k+\frac{1}{2}l_{ki}h_{ki}^j.
\end{equation}
Since $h_{ij}^k, h_{ki}^j$ are the signed distances, thus $A_i^{jk}$ is an algebraic sum of the area of triangles, i.e. a signed area.

\begin{figure}[!ht]
\centering
\begin{overpic}[scale=0.6]{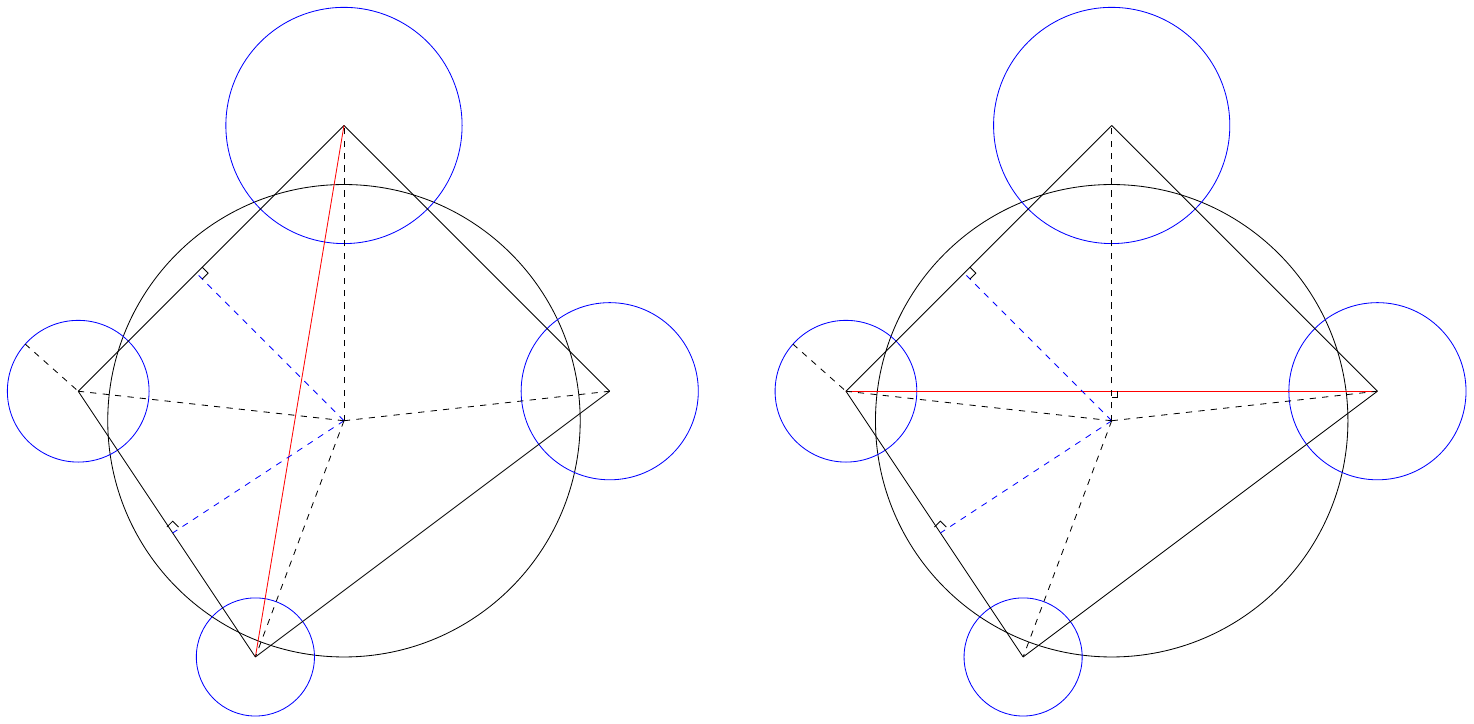}
\put (4,20) {$i$}
\put (41.5,20) {$l$}
\put (23,41) {$k$}
\put (17,2) {$j$}
\put (3.5,24.5) {$r_i$}
\put (14,17) {$h_{ij}^k$}
\put (16.5,27) {$h_{ki}^j$}
\put (23.5,19) {$c_{ijk}$}

\put (56.5,20) {$i$}
\put (94,20) {$l$}
\put (75,41) {$k$}
\put (69,2) {$j$}
\put (56,24.5) {$r_i$}
\put (76,19) {$c_{ijk}$}
\put (66.5,17) {$h_{ij}^l$}
\put (69,27) {$h_{ki}^l$}
\put (76,21) {$h_{il}^j$}
\end{overpic}
\caption{Weighted Delaunay triangulation $\mathcal{T}_1$ (left) and weighted Delaunay triangulation $\mathcal{T}_2$ (right).}
\label{figure 3}
\end{figure}
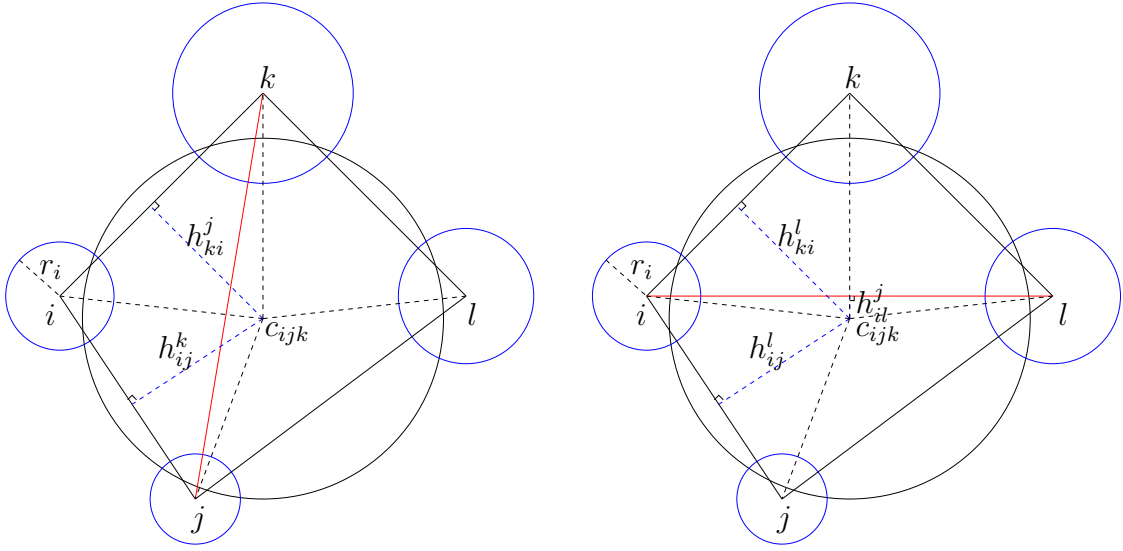

\begin{lemma}\label{Lem: independent}
The area $A_i$ is independent of the choice of weighted Delaunay triangulations of a decorated PE surface.
\end{lemma}
\proof
Suppose a decorated quadrilateral $\{ijlk\}$ is in a face of the weighted Delaunay tessellation of a decorated PE surface,
then there exist two weighted Delaunay triangulations $\mathcal{T}_1$ and $\mathcal{T}_2$ of the decorated PE surface
with an edge $\{jk\}$ in $\mathcal{T}_1$ flipped to another edge $\{il\}$ in $\mathcal{T}_2$. Please refer to Figure \ref{figure 3}.
We just need to prove the signed area $A^{jk}_i$ in $\mathcal{T}_1$ is equal to the signed area $A_i^{kl}+A_i^{jl}$ in $\mathcal{T}_2$.
In $\mathcal{T}_1$, the signed area at the vertex $i$ in $\{ijlk\}$ is
$A_i^{jk}=\frac{1}{2}l_{ki}h_{ki}^j+\frac{1}{2}l_{ij}h_{ij}^k$.
In $\mathcal{T}_2$,
the signed area at the vertex $i$ in $\{ijlk\}$ is
\begin{equation*}
\begin{aligned}
A_i^{kl}+A_i^{jl}
&=\frac{1}{2}l_{ki}h_{ki}^l+\frac{1}{2}l_{il}h_{il}^k
+\frac{1}{2}l_{ij}h_{ij}^l+\frac{1}{2}l_{il}h_{il}^j\\
&=\frac{1}{2}l_{ki}h_{ki}^l+\frac{1}{2}l_{ij}h_{ij}^l
+\frac{1}{2}l_{il}(h_{il}^k+h_{il}^j).
\end{aligned}
\end{equation*}
Since $\mathcal{T}_1$ and $\mathcal{T}_2$ are two weighted Delaunay triangulations of the same decorated PE metric on $(S,V)$, then $h_{il}^k+h_{il}^j=0$ by (\ref{Eq: F9}).
One can also refer to \cite{BL} (Proposition 3.4) for this.
Moreover, $h_{ki}^l=h_{ki}^j$ and $h_{ij}^l=h_{ij}^k$.
Then $A_i^{kl}+A_i^{jl}=A_i^{jk}$.
\qed

Denote $c_{ij}$ as the center of the edge $\{ij\}$, which is obtained by projecting the center $c_{ijk}$ to the line determined by $\{ij\}$.
Denote $d_{ij}$ as the signed distance of $c_{ij}$ to the vertex $i$,
which is positive if $c_{ij}$ is on the same side of $i$ as $j$ along the line determined by $\{ij\}$ and negative otherwise (or zero if $c_{ij}$ is the same as $i$).
In general, $d_{ij}\neq d_{ji}$.
Since the face-circle $C_{ijk}$ is orthogonal to the vertex-circle at the vertex $j$,
we have
\begin{equation}\label{Eq: F13}
r_{ijk}^2+r_j^2=d_{jk}^2+(h^i_{jk})^2=d_{ji}^2+(h^k_{ij})^2.
\end{equation}
Please refer to Figure \ref{figure 2} for this.
Moreover, we have the following explicit expressions of $d_{ij}$ and $h_{ij}^k$ due to Glickenstein \cite{Glickenstein JDG}, i.e.,
\begin{equation}\label{Eq: d ij}
d_{ij}=\frac{r_i^2+r_ir_jI_{ij}}{l_{ij}},
\end{equation}
and
\begin{equation}\label{Eq: h ijk}
h_{ij}^k=\frac{d_{ik}-d_{ij}\cos \theta_{jk}^i}{\sin \theta_{jk}^i},
\end{equation}
where $\theta^i_{jk}$ is the inner angle of the triangle $\{ijk\}$ at the vertex $i$.
The equality (\ref{Eq: d ij}) implies that $d_{ij}\in \mathbb{R}$ is independent of the existence of the center $c_{ijk}$.
Since each pair of vertex-circles is required to be separated, then $I>1$.
This implies
\begin{equation*}
d_{rs}>0,\ \ \forall\{r,s\}\subseteq\{i,j,k\}.
\end{equation*}
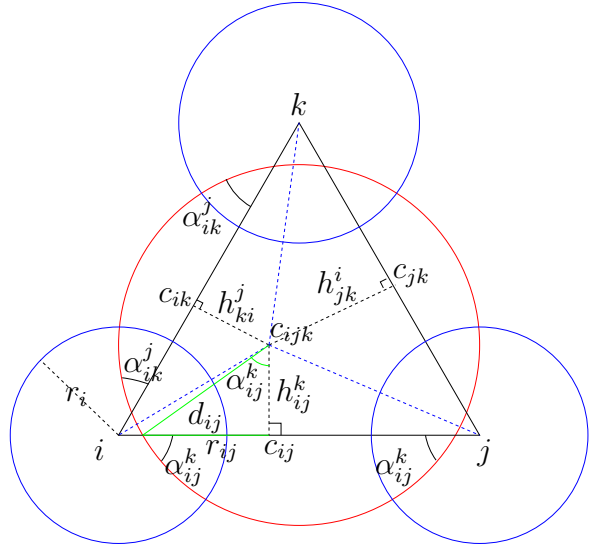
\begin{figure}[htbp]
\centering
\begin{overpic}[scale=0.4]{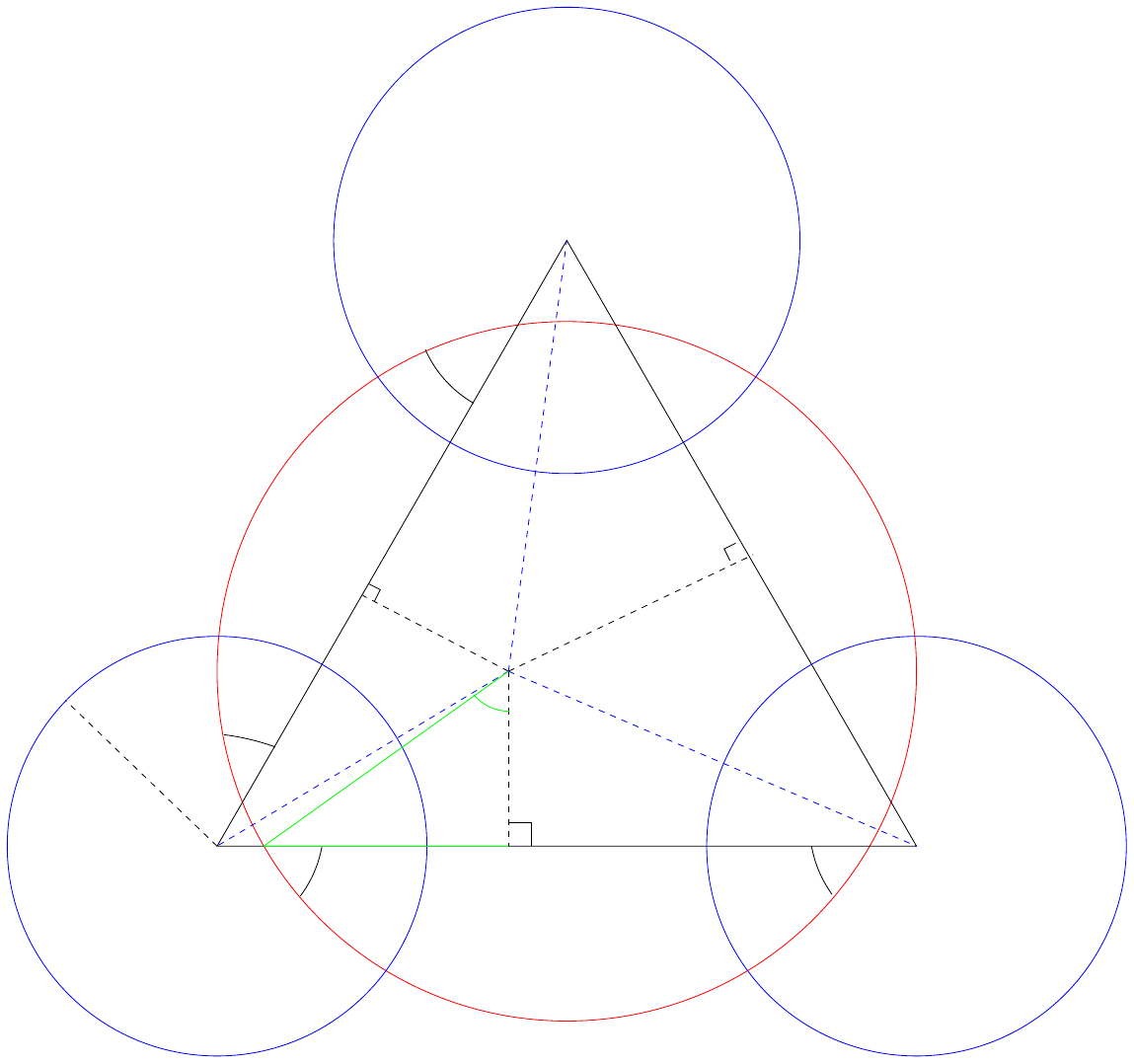}
\put (15,15) {$i$}
\put (80.5,15) {$j$}
\put (48.5,74) {$k$}
\put (30,55) {$\alpha_{ik}^j$}
\put (20,30) {$\alpha_{ik}^j$}
\put (37.5,27) {$\alpha_{ij}^k$}
\put (27,12.5) {$\alpha_{ij}^k$}
\put (63,12.5) {$\alpha_{ij}^k$}
\put (10,25) {$r_i$}
\put (46,25) {$h_{ij}^k$}
\put (36,40) {$h_{ki}^j$}
\put (53,43.5) {$h_{jk}^i$}
\put (31,21) {$d_{ij}$}
\put (45,36) {$c_{ijk}$}
\put (44,16) {$c_{ij}$}
\put (34,16) {$r_{ij}$}
\put (26,42) {$c_{ik}$}
\put (66,46) {$c_{jk}$}
\end{overpic}
\caption{Dates for a decorated triangle $\{ijk\}\in F$}
\label{figure 2}
\end{figure}

The following lemma gives some useful formulas.
\begin{lemma}[\cite{Guo,Xu AIM,Xu MRL}]
\label{Lem: useful 1}
Let $\{ijk\}$ be a decorated triangle with the edge lengths $l_{ij}, l_{jk}, l_{ki}$ defined by (\ref{Eq: DCE2}).
If the decorated triangle $\{ijk\}$ is non-degenerate, then
\begin{equation}\label{Eq: angle deform}
\frac{\partial \theta_{jk}^i}{\partial u_j}
=\frac{\partial \theta_{ki}^j}{\partial u_i}
=\frac{h_{ij}^k}{l_{ij}}, \ \ \
\frac{\partial \theta_{jk}^i}{\partial u_i}
=-\frac{\partial \theta_{jk}^i}{\partial u_j}
-\frac{\partial \theta_{jk}^i}{\partial u_k},
\end{equation}
where
\begin{equation}\label{Eq: h_ijk}
h_{ij}^k=\frac{r_i^2r_j^2r_k^2}{2A_{ijk}l_{ij}}
[\kappa_k^2(1-I_k^2)+\kappa_j\kappa_k\gamma_{i}
+\kappa_i\kappa_k\gamma_{j}]				=\frac{r_i^2r_j^2r_k^2}{2A_{ijk}l_{ij}}\kappa_kh_k			
\end{equation}
with $A_{ijk}=\frac{1}{2}l_{ij}l_{jk}\sin \theta_{ki}^j$, $\gamma_i=I_{jk}+I_{ij}I_{ki}$,\ $\kappa_i:=r_i^{-1}$ and
\begin{equation}\label{Eq: h_i}
\begin{aligned}
h_i=\kappa_i(1-I_{jk}^2)+\kappa_j\gamma_{k}+\kappa_k\gamma_{j}.
\end{aligned}
\end{equation}

\end{lemma}

As a direct application of Lemma \ref{Lem: useful 1}, we have the following result.
\begin{lemma}\label{Lem: A_ijk u_i}
The area $A_{ijk}(u)$ of each decorated triangle $\{ijk\}\in F$ is an analytic function with
\begin{equation}\label{Eq: A_ijk u_i}
\frac{\partial A_{ijk}}{\partial u_i}=A_i^{jk}.
\end{equation}
\end{lemma}
\proof
By (\ref{Eq: h ijk}), we have
\begin{equation*}
h^k_{ij}=\frac{d_{ik}-d_{ij}\cos \theta_{jk}^i}{\sin \theta_{jk}^i}, \quad
h_{ki}^j=\frac{d_{ij}-d_{ik}\cos \theta_{jk}^i}{\sin \theta_{jk}^i}.
\end{equation*}
Direct calculations give
\begin{equation}\label{Eq: F11}
h_{ki}^j=d_{ij}\sin \theta_{jk}^i-h^k_{ij}\cos \theta_{jk}^i.
\end{equation}
Combining (\ref{Eq: DCE2}), (\ref{Eq: inversive distance}) and (\ref{Eq: d ij}), it is easy to check that
\begin{equation}\label{Eq: F8}
\frac{\partial l_{ij}}{\partial u_i}=d_{ij}.
\end{equation}
Differentiating $A_{ijk}=\frac{1}{2}l_{ij}l_{jk}\sin \theta_{ki}^j$ with respect to $u_i$ gives
\begin{equation*}
\begin{aligned}
\frac{\partial A_{ijk}}{\partial u_i}
&=\frac{1}{2}\frac{\partial l_{ij}}{\partial u_i}l_{jk}\sin \theta_{ki}^j
+\frac{1}{2}l_{ij}l_{jk}\cos \theta_{ki}^j \frac{\partial \theta_{ki}^j}{\partial u_i}\\
&=\frac{1}{2}d_{ij}l_{jk}\sin \theta_{ki}^j
+\frac{1}{2}l_{ij}l_{jk}\cos \theta_{ki}^j \frac{h^k_{ij}}{l_{ij}}\\
&=\frac{1}{2}d_{ij}l_{ki}\sin \theta_{jk}^i
+\frac{1}{2}l_{jk}\cos \theta_{ki}^j h^k_{ij}\\
&=\frac{1}{2}d_{ij}l_{ki}\sin \theta_{jk}^i
+\frac{1}{2}(l_{ij}-l_{ki}\cos \theta_{jk}^i)h^k_{ij}\\
&=\frac{1}{2}l_{ki}(d_{ij}\sin \theta_{jk}^i-h^k_{ij}\cos \theta_{jk}^i)+\frac{1}{2}l_{ij}h^k_{ij}\\
&=\frac{1}{2}l_{ki}h_{ki}^j+\frac{1}{2}l_{ij}h_{ij}^k\\
&=A_i^{jk},
\end{aligned}
\end{equation*}
where the second equality uses (\ref{Eq: F8}) and (\ref{Eq: angle deform}), the third equality uses the sine laws and the penultimate line uses (\ref{Eq: F11}).
\qed

\begin{remark}
One can refer to Glickenstein \cite{Glickenstein JDG} for a nice geometric explanation of the result in Lemma \ref{Lem: A_ijk u_i}.
\end{remark}

\subsection{The extended energy function and the extended area function}
There exists a geometric relationship between the decorated triangle $\{ijk\}$ and the geometry of hyperbolic polyhedra in $3$-dimensional hyperbolic space.
Specially, there is a generalized hyperbolic tetrahedra in $\mathbb{H}^3$ with one ideal vertex and three hyper-ideal vertices corresponding to a decorated triangle $\{ijk\}$.
Please refer to \cite{BL} for more details on this fact.
Springborn \cite{Sp1} found the following explicit formula for the truncated volume $\mathrm{Vol}(ijk)$ of this generalized hyperbolic tetrahedra
\begin{equation}\label{Eq: volume}
\begin{aligned}
2\mathrm{Vol}(ijk)
=&\mathbb{L}(\theta_{jk}^i)+\mathbb{L}(\theta_{ki}^j)
+\mathbb{L}(\theta_{ij}^k)\\
&+\mathbb{L}(\frac{\pi+\alpha_{ki}^j+\alpha_{ij}^k-\theta_{jk}^i}{2})
+\mathbb{L}(\frac{\pi+\alpha_{ki}^j-\alpha_{ij}^k-\theta_{jk}^i}{2})\\
&+\mathbb{L}(\frac{\pi-\alpha_{ki}^j+\alpha_{ij}^k-\theta_{jk}^i}{2})
+\mathbb{L}(\frac{\pi-\alpha_{ki}^j-\alpha_{ij}^k-\theta_{jk}^i}{2})\\
&+\mathbb{L}(\frac{\pi+\alpha_{jk}^i+\alpha_{ij}^k-\theta_{ki}^j}{2})
+\mathbb{L}(\frac{\pi+\alpha_{jk}^i-\alpha_{ij}^k-\theta_{ki}^j}{2})\\
&+\mathbb{L}(\frac{\pi-\alpha_{jk}^i+\alpha_{ij}^k-\theta_{ki}^j}{2})
+\mathbb{L}(\frac{\pi-\alpha_{jk}^i-\alpha_{ij}^k-\theta_{ki}^j}{2})\\
&+\mathbb{L}(\frac{\pi+\alpha_{jk}^i+\alpha_{ki}^j-\theta_{ij}^k}{2})
+\mathbb{L}(\frac{\pi+\alpha_{jk}^i-\alpha_{ki}^j-\theta_{ij}^k}{2})\\
&+\mathbb{L}(\frac{\pi-\alpha_{jk}^i+\alpha_{ki}^j-\theta_{ij}^k}{2})
+\mathbb{L}(\frac{\pi-\alpha_{jk}^i-\alpha_{ki}^j-\theta_{ij}^k}{2}),
\end{aligned}
\end{equation}
where
\begin{equation}\label{Eq: Lobachevsky function}
\mathbb{L}(x)=-\int_0^x \log|2\sin(t)|dt
\end{equation}
is Milnor's Lobachevsky function.
Milnor's Lobachevsky function is bounded, odd, $\pi$-periodic and smooth except at integer multiples of $\pi$.
Please refer to \cite{Milnor, Rat} for more information on Milnor's Lobachevsky function $\mathbb{L}(x)$.

Set
\begin{equation}\label{Eq: F ijk}
\begin{aligned}
F_{ijk}(u_i,u_j,u_k)
=&-2\mathrm{Vol}(ijk)+\theta_{jk}^iu_i+\theta_{ki}^ju_j+\theta_{ij}^ku_k\\
&+(\frac{\pi}{2}-\alpha_{ij}^k)\lambda_{ij}
+(\frac{\pi}{2}-\alpha_{ki}^j)\lambda_{ki}
+(\frac{\pi}{2}-\alpha_{jk}^i)\lambda_{jk},
\end{aligned}
\end{equation}
where $\cosh \lambda_{ij}=I_{ij}$.
Then $\nabla F_{ijk}=(\theta_{jk}^i,\theta_{ki}^j, \theta_{ij}^k)$
and
\begin{equation}\label{Eq: property of F ijk}
F_{ijk}((u_i,u_j,u_k)+c(1,1,1))
=F_{ijk}(u_i,u_j,u_k)+c\pi
\end{equation}
for $c\in \mathbb{R}$.
Furthermore, on a decorated PE surface $(S,V,l,r)$ with a weighted Delaunay triangulation $\mathcal{T}$,
Bobenko-Lutz \cite{BL} defined the following function
\begin{equation}\label{Eq: F1}
\mathcal{H}_{\mathcal{T}}(u)
=\sum_{\{ijk\}\in F}F_{ijk}(u_i,u_j,u_k)
=-2\mathrm{Vol}(P_h)+\sum_{i\in V}\theta_iu_i+\sum_{\{ij\}\in E_{\mathcal{T}}}(\pi-\alpha_{ij})\lambda_{ij},
\end{equation}
where $P_h$ is the convex polyhedral cusp defined by the heights $h\in \mathbb{R}^V$, $\theta_i=\sum_{\{ijk\}\in F_\mathcal{T}}\theta^i_{jk}$ and $\alpha_{ij}=\alpha_{ij}^k+\alpha_{ij}^l$.
Note that the function $\mathcal{H}_{\mathcal{T}}(u)$ defined by (\ref{Eq: F1}) differs from its original definition in \cite{BL} (Equation 4-9) by some constant.
By (\ref{Eq: property of F ijk}), for $c\in \mathbb{R}$, we have
\begin{equation*}\label{Eq: property of H ijk}
\mathcal{H}_{\mathcal{T}}(u+c\mathbf{1})
=\mathcal{H}_{\mathcal{T}}(u)+c|F|\pi.
\end{equation*}
Using the function $\mathcal{H}_{\mathcal{T}}$, we define the following energy function
\begin{equation*}
\mathcal{E}_{\mathcal{T}}(u)
=-\mathcal{H}_{\mathcal{T}}(u)+2\pi\sum_{i\in V}u_i,
\end{equation*}
which is well-defined on $\mathcal{C}_\mathcal{T}(dist_{S},r)$ with
$\nabla_{u_i} \mathcal{E}_{\mathcal{T}}
=2\pi-\sum_{\{ijk\}\in F_\mathcal{T}}\theta^i_{jk}=W_i$.
Moreover, for $c\in \mathbb{R}$, we have
\begin{equation}\label{Eq: property of E}
\begin{aligned}
\mathcal{E}_{\mathcal{T}}(u+c\mathbf{1})
=&-\mathcal{H}_{\mathcal{T}}(u+c\mathbf{1})
+2\pi\sum_{i\in V}(u_i+c) \\
=&-\mathcal{H}_{\mathcal{T}}(u)-c|F|\pi
+2\pi\sum_{i\in V}u_i+2c|V|\pi \\
=&\mathcal{E}_{\mathcal{T}}(u)+2c\pi\chi(S),
\end{aligned}
\end{equation}
where $2|V|-|F|=2\chi(S)$ is used in the last line.

\begin{theorem}[\cite{BL}, Proposition 4.13]
\label{Thm: extended E}
For a discrete conformal factor $u\in \mathbb{R}^V$, let $\mathcal{T}$
be a weighted Delaunay triangulation of the decorated PE surface $(S,V,dist_S(u),r(u))$.
The map
\begin{equation}\label{Eq: extended H}
\begin{aligned}
\mathcal{H} :\  \mathbb{R}^V&\rightarrow \mathbb{R},\\
u&\mapsto \mathcal{H}_{\mathcal{T}}(u)
\end{aligned}
\end{equation}
is well-defined, concave, and twice continuously differentiable over $\mathbb{R}^V$.
\end{theorem}
Therefore, the function $\mathcal{E}_{\mathcal{T}}(u)$ defined on $\mathcal{C}_\mathcal{T}(dist_{S},r)$ can be extended to be
\begin{equation}\label{Eq: extended E}
\mathcal{E}(u)
=-\mathcal{H}(u)+2\pi\sum_{i\in V}u_i
=-\sum_{\{ijk\}\in F}F_{ijk}(u_i,u_j,u_k)
+2\pi\sum_{i\in V}u_i
\end{equation}
defined on $\mathbb{R}^V$.

\begin{definition}
Suppose $(S,V,\mathcal{T})$ is a triangulated surface with a decorated PE metric $(l,r)$.
The area function $A^\mathcal{T}_{tot}$ on $(S,V,\mathcal{T})$ is defined to be
\begin{equation*}
A^\mathcal{T}_{tot} :\  \mathcal{C}_\mathcal{T}(dist_{S},r)\rightarrow \mathbb{R},\\
\quad \ A^\mathcal{T}_{tot}(u)=\sum_{\{ijk\}\in F}A_{ijk}(u).
\end{equation*}
\end{definition}

By Lemma \ref{Lem: A_ijk u_i}, we have the following result.
\begin{corollary}\label{Cor: A_tot u_i}
The function $A^\mathcal{T}_{tot}$ is an analytic function with
\begin{equation}\label{Eq: A_tot u_i}
\frac{\partial A^\mathcal{T}_{tot}}{\partial u_i}=A_i.
\end{equation}
\end{corollary}

Lemma \ref{Lem: independent} and Corollary \ref{Cor: A_tot u_i} imply the following result, which shows the function $A^\mathcal{T}_{tot}$ defined on $\mathcal{C}_\mathcal{T}(dist_{S},r)$ can be extended.

\begin{theorem}\label{Thm: extended area}
For a discrete conformal factor $u\in \mathbb{R}^V$, let $\mathcal{T}$ be a weighted Delaunay triangulation of the decorated PE surface $(S,V,dist_S(u),r(u))$.
 The map
\begin{equation}\label{Eq: extended A_tot}
\begin{aligned}
A_{tot} :\  \mathbb{R}^V&\rightarrow \mathbb{R},\\
u&\mapsto A^\mathcal{T}_{tot}(u)
\end{aligned}
\end{equation}
is well-defined and once differentiable.
\end{theorem}
\proof
By Corollary \ref{Cor: A_tot u_i},
the function $A_{tot}$ is once differentiable in the interior of any $\mathcal{C}_\mathcal{T}(dist_{S},r)$.
At the boundary of $\mathcal{C}_\mathcal{T}(dist_{S},r)$,
the weighted triangulations induce the same weighted Delaunay tessellation.
The conclusion follows from Lemma \ref{Lem: independent}.
\qed

\section{The proof of Theorem \ref{Thm: main}}
\label{Sec: proof the main theorem}

\subsection{Variational principles with constraints}
In this subsection, we translate Theorem \ref{Thm: main} into an optimization problem with inequality constraints by variational principles,
which involves the function $\mathcal{E}$ defined by (\ref{Eq: extended E}).

\begin{proposition}\label{Prop: set A}
The set
\begin{equation*}
\mathcal{A}=\{u\in \mathbb{R}^V|A_{tot}(u)\leq 1\}.
\end{equation*}
is an unbounded closed subset of $\mathbb{R}^V$.
\end{proposition}
\proof
By Theorem \ref{Thm: extended area}, the set $\mathcal{A}$ is a closed subset of $\mathbb{R}^V$.
Since
$A_{ijk}((u_i,u_j,u_k)+(c,c,c))=e^{2c}A_{ijk}(u_i,u_j,u_k)$,
thus $A_{tot}(u+c\mathbf{1})=e^{2c}A_{tot}(u)$.
Then
$A_{tot}(u+c\mathbf{1})=e^{2c}A_{tot}(u)\leq 1$
is equivalent to
$c\leq -\frac{1}{2}\log A_{tot}(u)$.
This implies that the ray $\{u+c\mathbf{1}|c\leq -\frac{1}{2}\log A_{tot}(u)\}$ stays in the set $\mathcal{A}$.
Hence the set $\mathcal{A}$ is unbounded.
\qed

According to Proposition \ref{Prop: set A}, we have following result.
\begin{lemma}\label{Lem: minimum lies at the boundary}
If $\chi(S)<0$ and the function $\mathcal{E}(u)$ attains a minimum in the set $\mathcal{A}$,
then the minimum value point of $\mathcal{E}(u)$ lies at the boundary of $\mathcal{A}$, i.e.,
\begin{equation*}
\partial \mathcal{A}
=\{u\in \mathbb{R}^V|A_{tot}(u)=1\}.
\end{equation*}
Furthermore, there exists a decorated PE metric with constant discrete Gaussian curvature $K$ in the discrete conformal class.
\end{lemma}
\proof
Suppose the function $\mathcal{E}(u)$ attains a minimum at $u\in \mathcal{A}$.
Taking $c_0=-\frac{1}{2}\log A_{tot}(u)$,
then $c_0\geq0$ by $A_{tot}(u)\leq 1$.
By the proof of Proposition \ref{Prop: set A}, $u+c_0\mathbf{1}\in \mathcal{A}$.
Hence, by the additive property of the function $\mathcal{E}$ in (\ref{Eq: property of E}),
we have
\begin{equation*}
\mathcal{E}(u)\leq \mathcal{E}(u+c_0\mathbf{1})
=\mathcal{E}(u)+2 c_0\pi\chi(S).
\end{equation*}
This implies $c_0\leq0$ by $\chi(S)<0$.
Then $c_0=0$ and $A_{tot}(u)=1$.
Therefore, the minimum value point of $\mathcal{E}(u)$ lies in the set $\partial\mathcal{A}
=\{u\in\mathbb{R}^V|A_{tot}(u)=1\}$.
The conclusion follows from the following claim.\\
\textbf{Claim :} Up to scaling, the decorated PE metrics with constant discrete Gaussian curvature $K$ in the discrete conformal class
are in one-to-one correspondence with the critical points of the function $\mathcal{E}(u)$ under the constraint $A_{tot}(u)=1$.

We use the method of Lagrange multipliers to prove this claim.
Set
\begin{equation*}
G(u,\mu)=\mathcal{E}(u)-\mu(A_{tot}(u)-1),
\end{equation*}
where $\mu\in \mathbb{R}$ is a Lagrange multiplier.
If $u$ is a critical point of the function $\mathcal{E}$ under the constraint $A_{tot}(u)=1$,
then by (\ref{Eq: A_tot u_i}) and the fact $\nabla_{u_i} \mathcal{E}=W_i$, we have
\begin{equation*}
0=\frac{\partial G(u,\mu)}{\partial u_i}
=\frac{\partial \mathcal{E}(u)}{\partial u_i}
-\mu\frac{\partial A_{tot}(u)}{\partial u_i}
=W_i-\mu A_i.
\end{equation*}
This implies
\begin{equation*}
W_i=\mu A_i.
\end{equation*}
Since the anger defect $W$ defined by (\ref{Eq: curvature W}) satisfies the following discrete Gauss-Bonnet formula
\begin{equation*}
\sum_{i\in V} W_i=2\pi \chi(S),
\end{equation*}
we have
\begin{equation*}
2\pi \chi(S)=\sum_{i\in V}W_i
=\mu \sum_{i\in V}A_i=\mu A_{tot}=\mu.
\end{equation*}
under the constraint
$A_{tot}(u)=1$.
Therefore, the discrete Gaussian curvature
$$K_i=\frac{W_i}{A_i}=2\pi \chi(S)$$
for any $i\in V$.
\qed

\subsection{Reduction to Theorem \ref{Thm: key}}

By Lemma \ref{Lem: minimum lies at the boundary},
we just need to prove that the function $\mathcal{E}(u)$ attains the minimum in the set $\mathcal{A}$.
Recall the following classical result from calculus.

\begin{theorem}\label{Thm: calculus}
Let $A\subseteq \mathbb{R}^m$ be a closed set and $f: A\rightarrow \mathbb{R}$ be a continuous function. If every unbounded sequence $\{u_n\}_{n\in \mathbb{N}}$ in $A$ has a subsequence $\{u_{n_k}\}_{k\in \mathbb{N}}$ such that
$\lim_{k\rightarrow +\infty} f(u_{n_k})=+\infty$,
then $f$ attains a minimum in $A$.
\end{theorem}
One can refer to \cite{Kourimska Thesis} (Section 4.1) for a proof of Theorem \ref{Thm: calculus}.
The majority of the conditions in Theorem \ref{Thm: calculus} is satisfied,
including the set $\mathcal{A}$ is a closed subset of $\mathbb{R}^V$ by Proposition \ref{Prop: set A} and the function $\mathcal{E}$ is continuous by Theorem \ref{Thm: extended E}.
To prove Theorem \ref{Thm: main},
we just need to prove the following theorem.

\begin{theorem}\label{Thm: key}
If $\chi(S)<0$ and $\{u_n\}_{n\in \mathbb{N}}$ is an unbounded sequence in $\mathcal{A}$,
then there exists a subsequence $\{u_{n_k}\}_{k\in \mathbb{N}}$ of $\{u_n\}_{n\in \mathbb{N}}$ such that
$\lim_{k\rightarrow +\infty} \mathcal{E}(u_{n_k})=+\infty$.
\end{theorem}

\subsection{Behaviour of sequences of discrete conformal factors}

Let $\{u_n\}_{n\in \mathbb{N}}$ be an unbounded sequence in $\mathbb{R}^V$.
Denote its coordinate sequence at $j\in V$ by $\{u_{j,n}\}_{n\in \mathbb{N}}$.
Motivated by \cite{Kourimska}, we call the sequence $\{u_n\}_{n\in \mathbb{N}}$ with the following properties as a ``good" sequence.
\begin{description}
\item[(1)] It lies in one cell $\mathcal{C}_\mathcal{T}(dist_{S},r)$ of $\mathbb{R}^V$;
\item[(2)] There exists a vertex $i^*\in V$ such that $u_{i^*,n}\leq u_{j,n}$ for all $j\in V$ and $n\in \mathbb{N}$;
\item[(3)] Each coordinate sequence $\{u_{j,n}\}_{n\in \mathbb{N}}$ either converges, diverges properly to $+\infty$, or diverges properly to $-\infty$;
\item[(4)] For any $j\in V$, the sequence $\{u_{j,n}-u_{i^*,n}\}_{n\in \mathbb{N}}$ either converges or diverges properly to $+\infty$.
\end{description}

By Lemma \ref{Lem: finite decomposition},
it is obvious that every sequence of discrete conformal factors in $\mathbb{R}^V$ possesses a ``good" subsequence.
Hence, the ``good" sequence could be chosen without loss of generality.

In the following arguments, we use the following notations
\begin{equation}\label{Eq: F2}
l^n_{ij}=\sqrt{r^2_{i,n}+r^2_{j,n}+2 I_{ij}r_{i,n}r_{j,n}},
\end{equation}
\begin{equation}\label{Eq: F3}
r_{i,n}=e^{u_{i,n}}r_i,
\end{equation}
\begin{equation}\label{Eq: F4}
(l^n_{ij})^2
=(e^{2u_{i,n}}-e^{u_{i,n}+u_{j,n}})r^2_i
+(e^{2u_{j,n}}-e^{u_{i,n}+u_{j,n}})r^2_j
+e^{{u_{i,n}+u_{j,n}}}l_{ij}^2.
\end{equation}
For a decorated triangle $\{ijk\}\in F$ in $(S,V,\mathcal{T})$,
set
\begin{equation}\label{Eq: C ijk}
\mathcal{C}_{ijk}=\{(u_i,u_j,u_k)\in \mathbb{R}^3|u\in \mathcal{C}_\mathcal{T}(dist_{S},r)\}.
\end{equation}
Let $(u_{i,n},u_{j,n},u_{k,n})_{n\in \mathbb{N}}$ be a coordinate sequence in $\mathcal{C}_{ijk}$.
Then the edge lengths $l^n_{ij},l^n_{jk},l^n_{ki}$ satisfy the triangle inequalities for all $n\in \mathbb{N}$.

\begin{lemma}\label{Lem: two infty one bounded}
There exists no sequence in $\mathcal{C}_{ijk}$
such that as $n\rightarrow +\infty$,
\begin{equation*}
u_{r,n}\rightarrow +\infty, \quad u_{s,n}\rightarrow +\infty,
\quad u_{t,n}\leq C,
\end{equation*}
where $\{r,s,t\}=\{i,j,k\}$ and $C$ is a constant.
\end{lemma}
\proof
Without loss of generality, we assume
$\lim u_{i,n}=+\infty$, $\lim u_{j,n}=+\infty$
and the sequence $u_{k,n}\leq C_1$.
The equality (\ref{Eq: F3}) implies
$\lim r_{i,n}=+\infty$, $\lim r_{j,n}=+\infty$
and the sequence $r_{k,n}\leq C_2$.
Here $C_1,C_2$ are constants.
By (\ref{Eq: F2}), we have
\begin{equation*}
\begin{aligned}
(l^n_{jk}+l^n_{ki})^2
=&r^2_{i,n}+r^2_{j,n}+2r^2_{k,n}+2 I_{jk}r_{j,n}r_{k,n}
+2 I_{ki}r_{k,n}r_{i,n}\\
&+2\sqrt{(r^2_{j,n}+r^2_{k,n}+2 I_{jk}r_{j,n}r_{k,n})(r^2_{k,n}+r^2_{i,n}+2 I_{ki}r_{k,n}r_{i,n})}.
\end{aligned}
\end{equation*}
Note that $I_{ij}>1$, then
\begin{equation*}
\lim\frac{r^2_{k,n}+I_{jk}r_{j,n}r_{k,n}
+I_{ki}r_{k,n}r_{i,n}+\sqrt{(r^2_{j,n}+r^2_{k,n}+2 I_{jk}r_{j,n}r_{k,n})(r^2_{k,n}+r^2_{i,n}+2 I_{ki}r_{k,n}r_{i,n})}}{I_{ij}r_{i,n}r_{j,n}}<1.
\end{equation*}
Therefore, there exists $n\in \mathbb{N}$ such that
$(l^n)^2_{ij}=r^2_{i,n}+r^2_{j,n}+2 I_{ij}r_{i,n}r_{j,n}>(l^n_{jk}+l^n_{ki})^2$, i.e.,
$l^n_{ij}>l^n_{jk}+l^n_{ki}$.
This contradicts the triangle inequality $l^n_{ij}< l^n_{jk}+l^n_{ki}$.
\qed

Combining Lemma \ref{Lem: two infty one bounded} and the connectivity of the triangulation $\mathcal{T}$, we have the following result.
\begin{corollary}\label{Cor: one infty two converge}
For a discrete conformal factor $u\in \mathbb{R}^V$, let $\mathcal{T}$ be a weighted Delaunay triangulation of the decorated PE surface $(S,V,dist_S(u),r(u))$.
For any decorated triangle $\{ijk\}\in F$ in $\mathcal{T}$,
at least two of the three sequences $(u_{i,n}-u_{i^*,n})_{n\in \mathbb{N}}$,
$(u_{j,n}-u_{i^*,n})_{n\in \mathbb{N}}$,
$(u_{k,n}-u_{i^*,n})_{n\in \mathbb{N}}$
converge.
\end{corollary}

To characterize the function $F_{ijk}(u_i,u_j,u_k)$ in (\ref{Eq: F ijk}),
we need the following lemmas.

\begin{lemma}\label{Lem: converge 1}
Assume that the sequence $(u_{i,n})_{n\in \mathbb{N}}$ diverges properly to $+\infty$ and the sequences $(u_{j,n})_{n\in \mathbb{N}}$ and $(u_{k,n})_{n\in \mathbb{N}}$ converge.
Then the sequence $(\theta^{i,n}_{jk})_{n\in \mathbb{N}}$ converges to zero.
Furthermore, if the sequences $(\theta^{j,n}_{ki})_{n\in \mathbb{N}}$ and
$(\theta^{k,n}_{ij})_{n\in \mathbb{N}}$ converge to non-zero constants, then
\begin{description}
\item[(1)]
the sequences $(h_{jk}^{i,n})_{n\in \mathbb{N}}$,
$(h_{ki}^{j,n})_{n\in \mathbb{N}}$ and
$(h_{ij}^{k,n})_{n\in \mathbb{N}}$ converge;
\item[(2)]
the sequences $(\alpha^{i,n}_{jk})_{n\in \mathbb{N}}$,
$(\alpha^{j,n}_{ki})_{n\in \mathbb{N}}$ and
$(\alpha^{k,n}_{ij})_{n\in \mathbb{N}}$ converge.
\end{description}
\end{lemma}
\proof
By the assumption, we have
$\lim r_{i,n}=+\infty$, $\lim r_{j,n}=c_1$
and $\lim r_{k,n}=c_2$,
where $c_1,c_2$ are positive constants.
The equality (\ref{Eq: F2}) implies
\begin{equation}\label{Eq: F5}
\lim \frac{l_{ij}^n}{r_{i,n}}=1,\
\lim \frac{l_{ki}^n}{r_{i,n}}=1,\
\lim l_{jk}^n=c_3,
\end{equation}
where $c_3$ is a positive constant.
By the cosine law, we have
\begin{equation*}
\lim \cos\theta^{i,n}_{jk}
=\lim\frac{-(l_{jk}^n)^2+(l_{ij}^n)^2+(l_{ki}^n)^2}
{2l_{ij}^nl_{ki}^n}
=1.
\end{equation*}
This implies $\lim\theta^{i,n}_{jk}=0$.

Suppose the sequences $(\theta^{j,n}_{ki})_{n\in \mathbb{N}}$ and
$(\theta^{k,n}_{ij})_{n\in \mathbb{N}}$ converge to non-zero constants.
Then
\begin{equation}\label{Eq: A_ijk}
\lim \frac{A_{ijk}^n}{r_{i,n}}
=\lim \frac{l^n_{ij}l^n_{jk}\sin\theta^{j,n}_{ki}}{2r_{i,n}}
=c_4
\end{equation}
for some constant $c_4>0$.

\noindent\textbf{(1):}
Since $\kappa_i=\frac{1}{r_i}$, then
$\lim \kappa_{i,n}=0$, $\lim \kappa_{j,n}=\frac{1}{c_1}$
and $\lim \kappa_{k,n}=\frac{1}{c_2}$.
By (\ref{Eq: h_i}), we have
\begin{align*}
\lim h_{i,n}
=&\lim(\kappa_{i,n}(1-I_{jk}^2)+\kappa_{j,n}\gamma_{k}
+\kappa_{k,n}\gamma_{j})=c_5>0, \\
\lim h_{j,n}
=&\lim(\kappa_{j,n}(1-I_{ki}^2)+\kappa_{i,n}\gamma_{k}
+\kappa_{k,n}\gamma_i)=c_6, \\
\lim h_{k,n}
=&\lim(\kappa_{k,n}(1-I_{ij}^2)+\kappa_{i,n}\gamma_j
+\kappa_{j,n}\gamma_i)=c_7,
\end{align*}
where $c_5,c_6,c_7$ are constants.
Note that $c_6,c_7$ may be non-positive.
The equalities (\ref{Eq: h_ijk}) and (\ref{Eq: A_ijk}) imply
\begin{align*}
\lim h_{jk}^{i,n}
=&\lim \frac{r_{i,n}^2r_{j,n}^2r_{k,n}^2}{2A^n_{ijk}l^n_{jk}}
\kappa_{i,n}h_{i,n}
=\frac{c_1^2c^2_2c_5}{2c_3c_4}>0,\\
\lim h_{ki}^{j,n}
=&\lim \frac{r_{i,n}^2r_{j,n}^2r_{k,n}^2}{2A^n_{ijk}l^n_{ki}}
\kappa_{j,n}h_{j,n}
=\frac{c_1c^2_2c_6}{2c_4},\\
\lim h_{ij}^{k,n}
=&\lim \frac{r_{i,n}^2r_{j,n}^2r_{k,n}^2}{2A^n_{ijk}l^n_{ij}}
\kappa_{k,n}h_{k,n}	
=\frac{c_1^2c_2c_7}{2c_4}.
\end{align*}
Hence the sequences $(h_{jk}^{i,n})_{n\in \mathbb{N}}$,
$(h_{ki}^{j,n})_{n\in \mathbb{N}}$ and
$(h_{ij}^{k,n})_{n\in \mathbb{N}}$ converge.

\noindent\textbf{(2):}
The equality (\ref{Eq: d ij}) implies
\begin{equation}\label{Eq: F12}
\lim d^n_{jk}=\lim \frac{r_{j,n}^2+r_{j,n}r_{k,n}I_{jk}}{l^n_{jk}}
=\frac{c_1^2+c_1c_2I_{jk}}{c_3}>0.
\end{equation}
By (\ref{Eq: F13}), we have
\begin{equation*}
\lim (r^n_{ijk})^2
=\lim[(d^n_{jk})^2+(h^{i,n}_{jk})^2-r_{j,n}^2]
=c_8.
\end{equation*}
where $c_8$ is a constant.
Note that $h^{i}_{jk}=r_{ijk}\cos \alpha_{jk}^{i}$.
Hence,
\begin{align*}
\lim\cos\alpha_{jk}^{i,n}
=&\lim \frac{h_{jk}^{i,n}}{r^n_{ijk}}
=\frac{c_1^2c^2_2c_5}{2c_3c_4\sqrt{c_8}}>0,\\
\lim\cos\alpha_{ki}^{j,n}
=&\lim \frac{h_{ki}^{j,n}}{r^n_{ijk}}
=\frac{c_1c^2_2c_6}{2c_4\sqrt{c_8}},\\
\lim\cos\alpha_{ij}^{k,n}
=&\lim \frac{h_{ij}^{k,n}}{r^n_{ijk}}
=\frac{c_1^2c_2c_7}{2c_4\sqrt{c_8}}.
\end{align*}
Then the sequences $(\alpha^{i,n}_{jk})_{n\in \mathbb{N}}$,
$(\alpha^{j,n}_{ki})_{n\in \mathbb{N}}$ and
$(\alpha^{k,n}_{ij})_{n\in \mathbb{N}}$ converge.
\qed

\begin{lemma}\label{Lem: converge 2}
Assume that the sequence $(u_{i,n})_{n\in \mathbb{N}}$ diverges properly to $+\infty$ and the sequences $(u_{j,n})_{n\in \mathbb{N}}$ and $(u_{k,n})_{n\in \mathbb{N}}$ converge.
If the sequence $(\theta^{j,n}_{ki})_{n\in \mathbb{N}}$ converge to zero,
then
\begin{equation*}
\lim h^{i,n}_{jk}=+\infty,\
\lim h^{j,n}_{ki}=+\infty,\
\lim h^{k,n}_{ij}=-\infty.
\end{equation*}
\end{lemma}
\proof
Lemma \ref{Lem: converge 1} shows that $\lim\theta^{i,n}_{jk}=0$,
thus
$\lim (\theta^{j,n}_{ki}+\theta^{k,n}_{ij})=\pi$.
Since $\lim \theta^{j,n}_{ki}=0$,
then $\lim \theta^{k,n}_{ij}=\pi$.
Then
\begin{equation}\label{Eq: F10}
\lim \frac{A_{ijk}^n}{r_{i,n}}
=\lim \frac{l^n_{ij}l^n_{jk}\sin\theta^{j,n}_{ki}}{r_{i,n}}
=0.
\end{equation}
By the proof of Lemma \ref{Lem: converge 1}, we have
\begin{equation*}
\lim h_{jk}^{i,n}
=\lim \frac{r_{i,n}^2r_{j,n}^2r_{k,n}^2}{2A^n_{ijk}l^n_{jk}}
\kappa_{i,n}h_{i,n}
=\lim \frac{r_{i,n}^2c_1^2c_2^2}{2A^n_{ijk}c_3}\cdot
\frac{1}{r_{i,n}}c_5
=+\infty,
\end{equation*}
where (\ref{Eq: F10}) is used and $c_1,c_2,c_3,c_5$ are positive constants.
Similar to (\ref{Eq: F12}), we have
\begin{align*}
\lim d^n_{ji}
=&\lim\frac{r_{j,n}^2+r_{i,n}r_{j,n}I_{ij}}{l^n_{ij}}=c_9, \\
\lim d^n_{ki}
=&\lim\frac{r_{k,n}^2+r_{i,n}r_{k,n}I_{ki}}{l^n_{ki}}=c_{10}.
\end{align*}
Here $c_9,c_{10}$ are positive constants.
By (\ref{Eq: F13}), we have
\begin{align*}
(r^n_{ijk})^2
&=(d^n_{jk})^2+(h^{i,n}_{jk})^2-r_{j,n}^2 \\
&=(d^n_{ji})^2+(h^{k,n}_{ij})^2-r_{j,n}^2 \\
&=(d^n_{ki})^2+(h^{j,n}_{ki})^2-r_{k,n}^2.
\end{align*}
This implies $\lim r^n_{ijk}=+\infty$, $\lim (h^{k,n}_{ij})^2=+\infty$
and $\lim (h^{j,n}_{ki})^2=+\infty$.
Therefore, we have the following four cases
\begin{description}
\item[$(i)$] $\lim h^{k,n}_{ij}=+\infty$
and $\lim h^{j,n}_{ki}=+\infty$;
\item[$(ii)$] $\lim h^{k,n}_{ij}=-\infty$
and $\lim h^{j,n}_{ki}=-\infty$;
\item[$(iii)$] $\lim h^{k,n}_{ij}=+\infty$
and $\lim h^{j,n}_{ki}=-\infty$;
\item[$(iv)$] $\lim h^{k,n}_{ij}=-\infty$
and $\lim h^{j,n}_{ki}=+\infty$.
\end{description}
For the case $(i)$, since $\lim h_{jk}^{i,n}>0$, $\lim h^{k,n}_{ij}>0$ and $\lim h^{j,n}_{ki}>0$.
This implies that the center $c_{ijk}$ of the face-circle $C_{ijk}$ lies in the interior of the triangle $\{ijk\}$ by the definition of $h_{jk}^i, h^{k}_{ij}, h^{j}_{ki}$.
However, in this case, $\lim h_{jk}^{i,n}, \lim h^{k,n}_{ij},
\lim h^{j,n}_{ki}$ are bounded.
This is a contradiction.
\begin{figure}[htbp]
\centering
\begin{overpic}[scale=0.8]{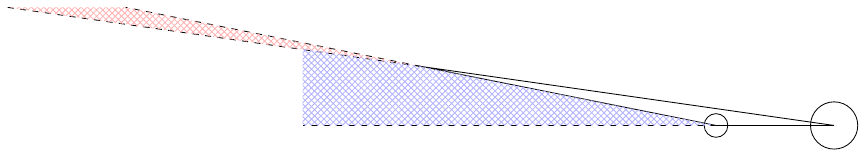}
\put(48,11){$i$}
\put(96,-2){$j$}
\put(82,-1){$k$}
\end{overpic}
\caption{Domain of the center $c_{ijk}$ in the decorated triangle on surface}
\label{figure 4}
\end{figure}
Both the cases $(ii)$ and $(iii)$ imply $d_{kj}<0$.
This contradicts with the fact that $d_{rs}>0$ for any $\{r, s\}\subseteq \{i,j,k\}$.
Indeed, the center $c_{ijk}$ lies in the red region in Figure \ref{figure 4} in the case $(ii)$ and  lies in the blue region in Figure \ref{figure 4} in the case $(iii)$.
By projecting the center $c_{ijk}$ to the line determined by $\{jk\}$, we have $d_{kj}<0$.
This completes the proof.
\qed

\begin{remark}\label{Rem: 2}
Similar to the proof of Lemma \ref{Lem: converge 2},
if the sequence $(\theta^{k,n}_{ij})_{n\in \mathbb{N}}$ converges to zero,
then
$\lim h^{i,n}_{jk}=+\infty,\
\lim h^{j,n}_{ki}=-\infty,\
\lim h^{k,n}_{ij}=+\infty.$
\end{remark}

\begin{figure}[htbp]
\centering
\begin{overpic}[scale=0.5]{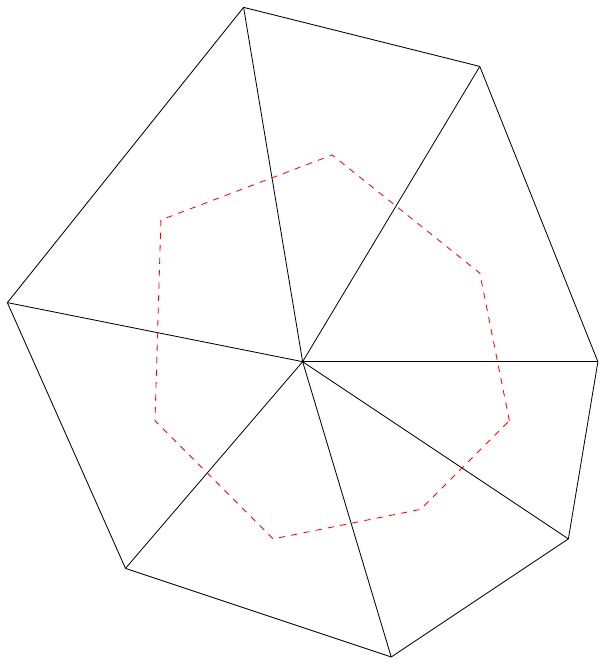}
\put(74,90){$j+1$}
\put(92,44){$j$}
\put(87,14){$j-1$}
\put(50,38){$i$}
\put(78,16){$0$}
\put(82,21){$\pi$}
\put(85,39){$0$}
\put(85,46){$\pi$}
\put(70,82){$0$}
\put(64,86){$\pi$}
\put(38,92){$0$}
\put(33,90){$\pi$}
\put(35,100){$s$}
\put(6,55){$0$}
\put(5,48){$\pi$}
\put(-4,50){$1$}
\put(18,18){$0$}
\put(23,14){$\pi$}
\put(16,8){$2$}
\put(53,4){$0$}
\put(59,5){$\pi$}
\end{overpic}
\caption{A star triangulation of a star-shaped $s$-sided polygon}
\label{figure 5}
\end{figure}
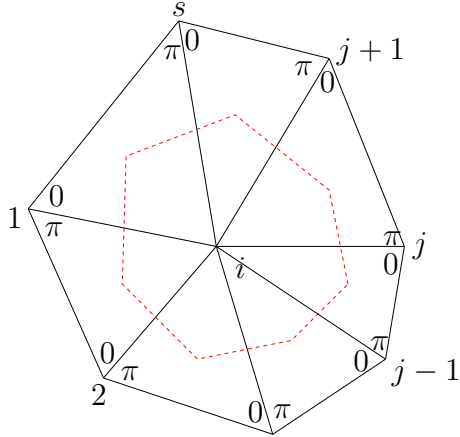

Consider a star-shaped $s$-sided polygon in the marked surface with boundary vertices $1,\cdots,s$ ordered cyclically ($v_{s+1}=v_1$).
Please refer to Figure \ref{figure 5}.
Let $i\in V$ be a vertex such that the sequence $(u_{i,n})_{n\in \mathbb{N}}$ diverges properly to $+\infty$ and the sequences $(u_{j,n})_{n\in \mathbb{N}}$ converge for $j\sim i$.

\begin{lemma}\label{Lem: converge 3}
The sequences of inner angles at the boundary vertices in the triangles of a star-shaped polygon converge to non-zero constants.
\end{lemma}
\proof
As $\lim u_{i,n}=+\infty$ and $\lim u_{j,n}=C$ for $j\sim i$.
By Lemma \ref{Lem: converge 1}, for any $j=1,...,s$,
we have $\lim \theta^{i,n}_{j-1,j}=0$ and hence
$\lim (\theta^{j-1,n}_{i,j}+\theta^{j,n}_{i,j-1})=\pi$.
We prove the result by contradiction.
Without loss of generality,
we assume $\lim \theta^{j-1,n}_{i,j}=\pi$
and $\lim \theta^{j,n}_{i,j-1}=0$ in the triangle $\{i,j-1,j\}$.
Then for $n$ large enough, we have
\begin{equation*}
l^n_{i,j-1}<l^n_{i,j}.
\end{equation*}
By Lemma \ref{Lem: converge 2},
we have $\lim h^{j,n}_{i,j-1}=+\infty,\
\lim h^{j-1,n}_{i,j}=-\infty,\
\lim h^{i,n}_{j-1,j}=+\infty$.
Since the edge $\{i,j\}$ is weighted Delaunay,
thus by (\ref{Eq: F9}), we have
\begin{equation*}
h^{j-1, n}_{i,j}+h^{j+1, n}_{i,j}\geq 0.
\end{equation*}
This implies $\lim h^{j+1,n}_{i,j}=+\infty$.

In the triangle $\{i,j,j+1\}$, suppose the sequences
$(\theta^{j,n}_{i,j+1})_{n\in \mathbb{N}}$ and
$(\theta^{j+1,n}_{i,j})_{n\in \mathbb{N}}$ converge to non-zero constants.
By Lemma \ref{Lem: converge 1},
the sequences $(h_{i,j}^{j+1,n})_{n\in \mathbb{N}}$ and
$(h_{i,j+1}^{j,n})_{n\in \mathbb{N}}$ converge.
This contradicts $\lim h^{j+1,n}_{i,j}=+\infty$.
Hence the sequences
$(\theta^{j,n}_{i,j+1})_{n\in \mathbb{N}}$ or
$(\theta^{j+1,n}_{i,j})_{n\in \mathbb{N}}$ converge to zero.
By Lemma \ref{Lem: converge 2} and Remark \ref{Rem: 2},
we have $\lim \theta^{j,n}_{i,j+1}=\pi$,
$\lim \theta^{j+1,n}_{i,j}=0$ and
$\lim h^{j,n}_{i,j+1}=-\infty$.
Then for $n$ large enough, we have
\begin{equation*}
l^n_{i,j}< l^n_{i,j+1}.
\end{equation*}
Please refer to Figure \ref{figure 5}.
By induction, for $n$ large enough, we have
\begin{equation*}
l^n_{i,j-1}< l^n_{i,j}< l^n_{i,j+1}< l^n_{i,j+2}<...< l^n_{i,j-1}.
\end{equation*}
This is a contradiction.
\qed

Combining (\ref{Eq: F ijk}), Lemma \ref{Lem: converge 1} and Lemma \ref{Lem: converge 3}, we have the following result.
\begin{corollary}\label{Cor: F converge}
Assume that the sequence $(u_{i,n})_{n\in \mathbb{N}}$ diverges properly to $+\infty$ and the sequences $(u_{j,n})_{n\in \mathbb{N}}$ and $(u_{k,n})_{n\in \mathbb{N}}$ converge.
Then the sequence $(F_{ijk}(u_{i,n},u_{j,n},u_{k,n}))_{n\in \mathbb{N}}$ converges.
\end{corollary}
\proof
By the definition of $F_{ijk}(u_i,u_j,u_k)$ in (\ref{Eq: F ijk}), we have
\begin{equation*}
\begin{aligned}
F_{ijk}(u_{i,n},u_{j,n},u_{k,n})
=&-2\mathrm{Vol}^n(ijk)+\theta_{jk}^{i,n}u_{i,n}
+\theta_{ki}^{j,n}u_{j,n}+\theta_{ij}^{k,n}u_{k,n}\\
&+(\frac{\pi}{2}-\alpha_{ij}^{k,n})\lambda_{ij}
+(\frac{\pi}{2}-\alpha_{ki}^{j,n})\lambda_{ki}
+(\frac{\pi}{2}-\alpha_{jk}^{i,n})\lambda_{jk}.
\end{aligned}
\end{equation*}
Combining Lemma \ref{Lem: converge 1} and Lemma \ref{Lem: converge 3} gives that
$\lim\theta^{i,n}_{jk}=0$, the sequences $(\theta^{j,n}_{ki})_{n\in \mathbb{N}}$ and
$(\theta^{k,n}_{ij})_{n\in \mathbb{N}}$ converge to non-zero constants and
the sequences $(\alpha^{i,n}_{jk})_{n\in \mathbb{N}}$,
$(\alpha^{j,n}_{ki})_{n\in \mathbb{N}}$ and
$(\alpha^{k,n}_{ij})_{n\in \mathbb{N}}$ converge.
Combining the continuity of Milnor's Lobachevsky function defined by (\ref{Eq: Lobachevsky function}) and the definition of the truncated volume $\mathrm{Vol}(ijk)$ defined by (\ref{Eq: volume}), we have that the sequence $(\mathrm{Vol}^n(ijk))_{n\in \mathbb{N}}$ converges.
Note that $\lambda_{ij}=\mathrm{arccosh} I_{ij}$ keeps invariant.
Hence,
$$\lim F_{ijk}(u_{i,n},u_{j,n},u_{k,n})=\lim \theta_{jk}^{i,n}u_{i,n}+c_{11}$$
for some constant $c_{11}$.
By (\ref{Eq: F3}), we have
$u_{i,n}=\log r_{i,n}-\log r_i$.
Then
\begin{equation*}
\begin{aligned}
\lim\theta_{jk}^{i,n}u_{i,n}
=&\lim \sin\theta^{i,n}_{jk}(\log r_{i,n}-\log r_i)\\
=&\lim\frac{2A_{ijk}^n}{l_{ij}^nl_{ki}^n}\log r_{i,n}\\
=&2c_4\lim\frac{\log r_{i,n}}{r_{i,n}}\\
=&0,
\end{aligned}
\end{equation*}
where the equalities (\ref{Eq: F5}) and (\ref{Eq: A_ijk}) is used in the second line and $\lim_{x\rightarrow +\infty}\frac{1}{x}\log x=0$ is used in the third line.
Therefore, $\lim F_{ijk}(u_{i,n},u_{j,n},u_{k,n})=c_{11}$.
\qed

The following lemma gives an asymptotic expression of the function $\mathcal{E}$.
\begin{lemma}\label{Lem: E decomposition}
There exists a convergent sequence $\{D_n\}_{n\in \mathbb{N}}$ such that the function $\mathcal{E}$ satisfies
\begin{equation*}
\mathcal{E}(u_n)=D_n+2\pi\left(u_{i^*,n}\chi(S)+\sum_{j\in V}(u_{j,n}-u_{i^*,n})\right).
\end{equation*}
\end{lemma}
\proof
By (\ref{Eq: extended E}), we have
\begin{align*}
\mathcal{E}(u_n)
=&-\sum_{\{ijk\}\in F}F_{ijk}(u_{i,n},u_{j,n},u_{k,n})
+2\pi\sum_{j\in V}u_{j,n}\\
=&-\sum_{\{ijk\}\in F}F_{ijk}
((u_{i,n},u_{j,n},u_{k,n})-u_{i^*,n}(1,1,1))
-\pi|F|u_{i^*,n}+2\pi\sum_{j\in V}u_{j,n}\\
=&D_n-\pi(2|V|-2\chi(S))u_{i^*,n}+2\pi\sum_{j\in V}u_{j,n}\\
=&D_n+2\pi\left(u_{i^*,n}\chi(S)+\sum_{j\in V}(u_{j,n}-u_{i^*,n})\right),
\end{align*}
where $D_n=-\sum_{\{ijk\}\in F}F_{ijk}
((u_{i,n},u_{j,n},u_{k,n})-u_{i^*,n}(1,1,1))$, the equation (\ref{Eq: property of F ijk}) is used in the second line and $2|V|-|F|=2\chi(S)$ is used in the third line.
The sequence $\{D_n\}_{n\in \mathbb{N}}$ converges by Corollary \ref{Cor: one infty two converge} and Corollary \ref{Cor: F converge}.
\qed

The following lemma gives the influence of the sequence $(u_n)_{n\in \mathbb{N}}$ on the area $A_{ijk}$ of a decorated triangle $\{ijk\}$.
\begin{lemma}\label{Lem: A ijk decomposition}
For a discrete conformal factor $u\in \mathbb{R}^V$, let $\mathcal{T}$ be a weighted Delaunay triangulation of the decorated PE surface $(S,V,dist_S(u),r(u))$.
Assume the sequences $(u_{j,n}-u_{i^*,n})_{n\in \mathbb{N}}$,
$(u_{k,n}-u_{i^*,n})_{n\in \mathbb{N}}$ converge for $\{ijk\}$  in  $\mathcal{T}$ with edge lengths $l_{ij}^n, l_{jk}^n, l_{ki}^n$.
\begin{description}
\item[(a)]
If $(u_{i,n}-u_{i^*,n})_{n\in \mathbb{N}}$ converges, there exists a convergent sequence of real numbers $(C_n)_{n\in \mathbb{N}}$ such that
\begin{equation}\label{Eq: key 1}
\log A_{ijk}^n=C_n+2u_{i^*,n}.
\end{equation}
\item[(b)]
If $(u_{i,n}-u_{i^*,n})_{n\in \mathbb{N}}$ diverges to $+\infty$, there exists a convergent sequence of real numbers $(C_n)_{n\in \mathbb{N}}$ such that
\begin{equation}\label{Eq: key 2}
\log A_{ijk}^n=C_n+u_{i,n}+u_{i^*,n}.
\end{equation}
\end{description}
\end{lemma}
\proof
Applying (\ref{Eq: F4}) to $A_{ijk}=\frac{1}{2}l_{ij}l_{jk}\sin \theta_{ki}^j$ gives
\begin{equation*}
\begin{aligned}
A^n_{ijk}
=&\frac{1}{2}l_{ij}^nl_{jk}^n\sin\theta_{ki}^{j,n}\\
=&\frac{1}{2}\sin\theta_{ki}^{j,n}
\sqrt{r^2_i e^{2u_{i,n}}+r^2_je^{2u_{j,n}}+(l_{ij}^2
-r^2_i-r^2_j)e^{(u_{i,n}+u_{j,n})}}\\
&\times \sqrt{r^2_j e^{2u_{j,n}}+r^2_ke^{2u_{k,n}}+(l_{jk}^2
-r^2_j-r^2_k)e^{(u_{j,n}+u_{k,n})}}.
\end{aligned}
\end{equation*}
Then
\begin{equation*}
\begin{aligned}
\log A^n_{ijk}
=&\log(\frac{1}{2}\sin\theta_{ki}^{j,n})+2u_{i^*,n}\\
&+\frac{1}{2}\log(r^2_i e^{2(u_{i,n}-u_{i^*,n})}
+r^2_je^{2(u_{j,n}-u_{i^*,n})}+(l_{ij}^2
-r^2_i-r^2_j)e^{(u_{i,n}-u_{i^*,n})+(u_{j,n}-u_{i^*,n})})\\
&+\frac{1}{2}\log(r^2_j e^{2(u_{j,n}-u_{i^*,n})}
+r^2_ke^{2(u_{k,n}-u_{i^*,n})}+(l_{jk}^2
-r^2_j-r^2_k)e^{(u_{j,n}-u_{i^*,n})+(u_{k,n}-u_{i^*,n})})\\
=&\log(\frac{1}{2}\sin\theta_{ki}^{j,n})+u_{i,n}+u_{i^*,n}\\
&+\frac{1}{2}\log(r^2_i+r^2_je^{2(u_{j,n}-u_{i^*,n})-2(u_{i,n}-u_{i^*,n})}
+(l_{ij}^2-r^2_i-r^2_j)
e^{-(u_{i,n}-u_{i^*,n})+(u_{j,n}-u_{i^*,n})})\\
&+\frac{1}{2}\log(r^2_j e^{2(u_{j,n}-u_{i^*,n})}
+r^2_ke^{2(u_{k,n}-u_{i^*,n})}+(l_{jk}^2
-r^2_j-r^2_k)e^{(u_{j,n}-u_{i^*,n})+(u_{k,n}-u_{i^*,n})}).
\end{aligned}
\end{equation*}
If the sequence $(u_{i,n}-u_{i^*,n})_{n\in \mathbb{N}}$ converges,
then $\log A_{ijk}^n=C_n+2u_{i^*,n}$.
If the sequence $(u_{i,n}-u_{i^*,n})_{n\in \mathbb{N}}$ diverges to $+\infty$,
then the sequence $(\theta_{ki}^{j,n})_{n\in \mathbb{N}}$ converges to a non-zero constant in $(0, \pi)$ by Lemma \ref{Lem: converge 3}.
This implies
$\log A_{ijk}^n
=C_n+u_{i,n}+u_{i^*,n}$.
In both cases, the sequence $(C_n)_{n\in \mathbb{N}}$ converges.
\qed

\subsection{Proof of Theorem \ref{Thm: key}}
Let $\{u_n\}_{n\in \mathbb{N}}$ be an unbounded ``good" sequence.
Suppose $\chi(S)<0$ and $\{u_n\}_{n\in \mathbb{N}}$ is an unbounded sequence in $\mathcal{A}$.
Combining $\chi(S)<0$ and Lemma \ref{Lem: E decomposition},
we just need to prove that $\lim_{n\rightarrow +\infty} u_{i^*,n}=-\infty$.
By the definition of ``good" sequence,
the sequence $\left(\sum_{j\in V}(u_{j,n}-u_{i^*,n})\right)_{n\in \mathbb{N}}$ converges to a finite number or diverges properly to $+\infty$.

If $\left(\sum_{j\in V}(u_{j,n}-u_{i^*,n})\right)_{n\in \mathbb{N}}$ converges to a finite number,
then the sequence $(u_{j,n}-u_{i^*,n})_{n\in \mathbb{N}}$ converges for all $j\in V$.
Since the sequence $\{u_n\}_{n\in \mathbb{N}}$ lies in $\mathcal{A}$, the area $A_{ijk}$ of each triangle is bounded from above.
This implies $\{u_{i^*,n}\}_{n\in \mathbb{N}}$ is bounded from above by (\ref{Eq: key 1}).
Then $\{u_{i^*,n}\}_{n\in \mathbb{N}}$ converges to a finite number or diverges properly to $-\infty$.
Suppose $\{u_{i^*,n}\}_{n\in \mathbb{N}}$ converges to a finite number.
Since $(u_{j,n}-u_{i^*,n})_{n\in \mathbb{N}}$ converges for all $j\in V$,
then $\{u_{j,n}\}_{n\in \mathbb{N}}$ are bounded for all $j\in V$, which implies $\{u_n\}_{n\in \mathbb{N}}$ is bounded.
This contradicts the assumption that $\{u_n\}_{n\in \mathbb{N}}$ is unbounded.
Therefore, the sequence $\{u_{i^*,n}\}_{n\in \mathbb{N}}$ diverges properly to $-\infty$.

If $\left(\sum_{j\in V}(u_{j,n}-u_{i^*,n})\right)_{n\in \mathbb{N}}$ diverges properly to $+\infty$,
then there exists at least one vertex $i\in V$ such that the sequence $(u_{i,n}-u_{i^*,n})_{n\in \mathbb{N}}$ diverges properly to $+\infty$.
By Corollary \ref{Cor: one infty two converge},
the sequences $(u_{j,n}-u_{i^*,n})_{n\in \mathbb{N}}$ and $(u_{k,n}-u_{i^*,n})_{n\in \mathbb{N}}$ converge for $j\sim i$ and $k\sim i$.
Since the area $A_{ijk}$ of each triangle is bounded from above,
thus $u_{i,n}+u_{i^*,n}\leq C$ and $u_{j,n}+u_{i^*,n}\leq C$ by (\ref{Eq: key 2}), where $C$ is a constant.
Then $(u_{i,n}-u_{i^*,n})+2u_{i^*,n}\leq C$.
This implies $\{u_{i^*,n}\}_{n\in \mathbb{N}}$ diverges properly to $-\infty$.
\qed

\begin{remark}\label{Rem: difficult}
For the case $\chi(S)>0$, Kou\v{r}imsk\'{a} \cite{Kourimska,Kourimska Thesis} gave the existence of PE metrics with constant discrete Gaussian curvatures.
However, we can not get similar results.
The main difference is that the edge length defined by  (\ref{Eq: DCE2}) involves the square term of discrete conformal factors, such as $e^{2u_i}$,
while the edge length defined by the vertex scalings only involves the mixed product of the first order terms, i.e., $e^{u_i+u_j}$.
Indeed, in this case, we can define the set $\mathcal{A}_+=\{u\in \mathbb{R}^V|A_{tot}(u)\geq 1\}$,
which is an unbounded closed subset of $\mathbb{R}^V$.
Under the conditions that $\chi(S)>0$ and the function $\mathcal{E}(u)$ attains a minimum in the set $\mathcal{A}_+$, Lemma \ref{Lem: minimum lies at the boundary} still holds.
Using Theorem \ref{Thm: calculus}, we just need to prove Theorem \ref{Thm: key} under the condition $\chi(S)>0$.
However, we can not get a good asymptotic expression of the area $A_{ijk}$.
The asymptotic expression of the area $A_{ijk}$ in (\ref{Eq: key 2}) involves $u_{i,n}+u_{i^*,n}$, which is not enough for this case.
\end{remark}

\end{document}